\documentclass[10pt]{article}
\usepackage[francais,english]{babel}
\usepackage{amsmath}
\usepackage{amsfonts}
\usepackage{amssymb}
\usepackage{amsthm}
\usepackage{graphics}
\usepackage{amscd}
\usepackage{fullpage}
\usepackage{epsfig}
\usepackage{color}
\newcommand{\Z}{\mathbb{Z}}
\newcommand{\R}{\mathbb{R}}
\newcommand{\N}{\mathbb{N}}

\newcommand{\C}{\mathbb{C}}

\newcommand{\U}{\mathbb{U}}
\newcommand{\E}{\mathbb{E}}

\newcommand{\mc}{\mathcal}

\newcommand{\eps}{\varepsilon}
\newcommand{\ind}{{\bf 1}}
\renewcommand{\P}{\mathbb{P}}

\renewcommand{\H}{\mathbb{H}}

\DeclareMathOperator{\Tr}{Tr}

\DeclareMathOperator{\SLE}{SLE}

\title{Euler integrals for commuting SLEs}
\author{Julien Dub\'edat\footnote{Courant Institute}}
\newtheorem{thm}{Theorem}

\newtheorem{Prop}[thm]{Proposition}
\newtheorem{Lem}[thm]{Lemma}

\begin{document}
\maketitle
\begin{abstract}
Schramm-Loewner Evolutions (SLEs) have proved an efficient way to
describe a single continuous random conformally invariant interface in
a simply-connected planar domain; the admissible probability
distributions are parameterized by a single positive parameter
$\kappa$. As shown in \cite{D6}, the coexistence of $n$ interfaces in
such a domain implies algebraic (``commutation'') conditions. In the
most interesting situations, the admissible laws on systems of $n$
interfaces are parameterized by $\kappa$ and the solution of a
particular (finite rank) holonomic system.

The study of solutions of differential systems, in particular their
global behaviour, often involves the use of integral
representations. In the present article, we provide Euler integral
representations for solutions of holonomic systems arising from SLE
commutation. Applications to critical percolation (general crossing
formulae), Loop-Erased Random Walks (direct derivation of Fomin's
formulae in the scaling limit), and Uniform Spanning Trees are
discussed. The connection with conformal restriction and Poissonized
non-intersection for chordal SLEs is also studied.
\end{abstract}

\section{Introduction}

Critical systems in the plane (such as critical percolation,
self-avoiding walk, ...) are generally conjectured to have scaling
limits satisfying certain conformal invariance properties. Schramm-Loewner
Evolutions (SLE), introduced in \cite{S0}, have proved a powerful tool
to describe and study these limits.

In several of these discrete models, such as percolation, one can
define more than one interface. Scaling limits of these systems of
interfaces should have similar conformal invariance properties. If one
studies the limit of a single simple path, then conformal invariance
and a ``Markov'' property specify the limit, if it exists, up to a
positive parameter, denoted by $\kappa$; this limit is $\SLE_\kappa$.

If one considers the joint law of several random curves, then the
limit can be encoded through Loewner's equations as (one-dimensional) diffusions. The
coefficients now depend on several parameters describing the boundary
conditions. But these coefficients have to satisfy some compatibility
conditions, expressed as commutation relations for associated
infinitesimal generators. These commutation conditions can be recast
as systems of linear PDEs with meromorphic coefficients and regular
singularities along hyperplanes (\cite{D6}). The key idea here is that
we want to define distributions on systems of geometric paths,
independently of any (artificial) time parameterization.

A particularly interesting set-up is the following: $(2n)$ points are
marked on the boundary of a planar simply-connected domain, and these
$(2n)$ points are joined by $n$ curves. Then it is shown in \cite{D6}
that the possible systems of random curves can be obtained by solving
a certain holonomic system, which is parameterized by $\kappa$ and
$n$. 

In this paper, we shall be  mainly interested in the properties of this family
of holonomic systems introduced in \cite{D6}. Solutions of this system can be related to particular
(local) martingales of SLE. These local martingales are essentially the Girsanov
densities for the systems of random curves w.r.t. independent
SLEs. Our main goal here is to
give integral representations of these solutions. When studying global
properties of solutions of differential equations
(asymptotics/boundary conditions, monodromy, \dots), use of such
integral representations is crucial. Also, although the theory of
holonomic systems is rich and deep, getting solutions of a particular
system is often challenging.

Cardy's formula for critical percolation (\cite{Ca,Sm1}) gives the
probability that there exists a crossing from left to right in a
rectangle. For $\kappa=6$, the system generalizes Cardy's formula to
situations were sides of a $(2n)$-gon are set to alternate colors.
Other related formulae have been proposed or proved, such as a formula
for expected number of crossings (\cite{Ca1}),  for crossing of
annuli (\cite{Ca2}), Watts' formula
for rectangles (\cite{Wa,D5}), and Pinson's formula for elliptic
curves (\cite{Pi}). These
formulae depend on one parameter (which is complex in Pinson's formula).

For $\kappa\in (0,8/3)$, the solutions of the system can be
interpreted in the continuous limit in terms of non-intersection of
$n$ independent $\SLE_\kappa$ and independent ``loop-soups''
(\cite{LSW3,LW,W3}) with appropriate intensity. If $\kappa=8/3$, this
loop-soup is empty, and we are studying the probability that $n$
independent $\SLE_{8/3}$ do not intersect. As $\SLE_{8/3}$ is
conjectured to be the scaling limit of Self-Avoiding Walks
(\cite{LSW5}), one can think in terms of the scaling limit of
non-intersecting self-avoiding walks.

In the case $\kappa=2$,  one can see
that the solutions of the system are exceptionally rational and
connect them to Fomin's determinantal formulae for Loop-Erased Random Walks
(see \cite{LSW2,Fomin,KozL}). This can also be seen as a particular
case of the previous restriction construction (in the continuum). In
the case $\kappa=8$, corresponding to the scaling limit of the
Uniform Spanning Tree (UST, see \cite{LSW2}), one can study the situation where
alternate wired and free conditions on the boundary of a domain force
the creation of $n$ non-intersecting Peano paths on the Manhattan
lattice. The scaling limit can be identified through Wilson's
algorithm.

It has been brought to our attention that the solutions obtained here
appear in the so-called Coulomb Gas formalism (see \cite{DotFat},
Chapter 9 in \cite{DiF}, and references therein), where they are
referred to as Coulomb Gas representations.
In this framework, deriving integral representations for correlation functions involves
(chiral) vertex operators and screening operators, and the
Feigin-Fuchs integral representation. Further developments along these
lines pertain to monodromy, fusion, and BRST cohomology.

In the present article, we obtain directly these integral representations
starting from the associated holonomic system. The analysis of this
system leads for example to results on rank and reductions using symmetry (as in
Section 4.4). Also, we focus on the interplay between the algebraic,
analytic and probabilistic aspects of the problem. For instance,
delicate analytic questions regarding boundary values can be bypassed
by probabilistic arguments (Section 6). In the case of Loop-Erased
Random Walks and Uniform Spanning Trees ($c=-2$), direct connections
with the combinatorial models are also established (3.2, 3.3).

The set-up, objects and notations are essentially the same as in
\cite{D6} and are mainly the standard SLE notations ($(g_t)$,
$\kappa$, \dots), so we will recall them only briefly. This article
can be read either as a self-contained study of a particular family of
holonomic systems, or as a companion paper to \cite{D6}.

The paper is organized as follows. After some background and
notations, several situations leading to natural systems of $n$ curves
are described. Discrete models, such as critical percolation,
loop-erased random walks (in relation with Fomin's formulae), and
uniform spanning trees (in relation with Wilson's algorithm) provide
important examples. In the continuous set-up, the theory of conformal
restriction and loop-soups (\cite{LSW3,LW}) can also be connected to
the problem of commuting SLEs (\cite{D6}). In the following section, as a preparation,
several particular cases (particular values of $\kappa$ and $n$) are
briefly studied; in particular, for $\kappa=6$ and $2n=6$, we get new
crossing probability expressions for critical percolation in conformal
hexagons. Building on these particular solutions, general (formal)
solutions are given, as Euler integral representations. We then
discuss the problem of identifying solutions corresponding to specific
geometric configurations; or, in other terms, how to choose the cycle
of integration in the Euler integral.

\section{Background and notations}

First we recall some definitions and fix notations. We shall be mainly
interested in chordal versions of SLE, i.e. random curves connecting
two boundary points of a plane simply-connected domain, satisfying
some conformal invariance properties.
For general background on SLE, see \cite{RS01,W1,Law}. Also, we will
use at some point results on the restriction property and the ``loop-soup'' (see \cite{LSW3,LW,W3}).

Consider the family of ODEs, indexed by $z$ in the upper half-plane
$\H=\{z:\Im z>0\}$:
$$\partial_tg_t(z)=\frac 2{g_t(z)-W_t}$$ 
with initial conditions $g_0(z)=z$, where $W_t$ is some real-valued
(continuous) function. These chordal Loewner equations are defined up
to explosion time $\tau_z$ (maybe infinite). Define:
$K_t=\overline{\{z\in\H:\tau_z<t\}}$.
Then $(K_t)_{t\geq 0}$ is an increasing family of compact subsets of
$\overline\H$; moreover, $g_t$ is the unique conformal equivalence
$\H\setminus K_t\rightarrow \H$ such that (hydrodynamic normalization
at $\infty$):
$$g_t(z)=z+o(1).$$ 
For any compact subset $K$ of $\overline{\H}$ such that $\H\setminus
K$ is simply-connected, we denote by $\phi$ the unique conformal
equivalence $\H\rightarrow \H\setminus K$ with hydrodynamic normalization
at $\infty$; so that $g_t=\phi_{K_t}$. The coefficient of $1/z$ in the Laurent expansion of $g_t$ at $\infty$
is by definition the half-plane capacity of $K_t$ at infinity; this
capacity equals $(2t)$.

If $W_t=x+\sqrt\kappa B_t$ where $(B_t)$ is a standard Brownian
Motion, then the Loewner chain $(K_t)$ (or the family $(g_t)$) defines
the chordal Schramm-Loewner Evolution with parameter $\kappa$ in
$(\H,x,\infty)$. The chain $K_t$ is generated by the trace $\gamma$, a
continuous process taking values in $\overline\H$, in the following
sense: $\H\setminus K_t$ is the unbounded connected component of
$\H\setminus\gamma_{[0,t]}$.

The trace is a continuous non self-traversing curve. It is a.s. simple
if $\kappa\leq 4$ and a.s. space-filling if $\kappa\geq 8$. Its
Hausdorff dimension is a.s. $1+\kappa/8$ if $\kappa\leq 8$ (and 2
otherwise).

If $(D,x,y)$ is a simply-connected domain with two marked points on the
boundary, $\SLE_\kappa$ from $x$ to $y$ in $D$ is defined as the image
of $\SLE_\kappa$ from 0 to $\infty$ in $\H$ (as defined above) by a
conformal equivalence $(\H,0,\infty)\rightarrow (D,x,y)$ (that exists
by Riemann's mapping theorem). With this definition, $\SLE$ satisfies
a ``Domain Markov'' property, which, together with conformal
equivalence, essentially characterizes it.

In \cite{D6}, the question of defining several $\SLE$ strands
simultaneously in a domain is addressed (say each of those
strands is absolutely continuous w.r.t. $\SLE_\kappa$). The key point
is that an appropriate ``Domain Markov'' condition for the joint
law imposes stark ``commutation conditions'' on the drift terms of the
driving processes. Elucidating those conditions is the main object of
\cite{D6}. We summarize the result in the following case: $(2n)$
points $x_1,\dots,x_{2n}$ are marked on the boundary of $\H$, and we
want to define jointly $n$ SLEs connecting these $(2n)$ points, such
that the joint distribution depends only on the position of the marked
points, is M\"obius invariant and satisfies the appropriate Markov
property. Then necessarily, the $\SLE$ growing at $x_i$ is driven by:
$$dX^{(i)}_t=\sqrt\kappa dB^{(i)}_t+\kappa\frac{\partial_i\varphi}\varphi(g_t(x_1),\dots,X^{(i)}_t,\dots,g_t(x_{2n}))dt$$
where $B^{(i)}$ is a standard Brownian motion and 
$\varphi$ is a non-vanishing function annihilated by the operators:

\begin{equation}\label{sys}
\left\{\begin{array}{l}\displaystyle\frac\kappa 2\partial_{kk}+\sum_{l\neq
k}\frac{2\partial_{l}}{x_l-x_k}+\frac{\kappa-6}\kappa\sum_{l\neq
k}\frac 1{(x_l-x_k)^2}, \hspace{1cm} k=1,\dots, 2n\\
\sum_k\partial_k\\
\sum_k x_k\partial_k-n(1-6/\kappa)\\
\sum_k x_k^2\partial_k-(1-6/\kappa)(x_1+\cdots+ x_{2n})
\end{array}\right.
\end{equation}

From a CFT point of view, this corresponds to differential equations
for correlation functions of $\phi_{1;2}$ primary fields.

If $n=1$, $\varphi(x_1,x_2)=(x_2-x_1)^{1-6/\kappa}$, and:
$$dX^{(1)}_t=\sqrt\kappa dB^{(1)}_t+\frac{\kappa-6}{X^{(1)}_t-g_t(x_2)}dt$$
which is ordinary chordal $\SLE_\kappa$ from $x_1$ to $x_2$ (with a
homographic change of coordinate).

The main question addressed in this article is how to obtain general
integral representations for solutions of this system and their
interpretation in terms of scaling limits for different discrete
models.

\section{Probabilistic interpretations}

In this section, we discuss probabilistic situations giving rise to
natural examples of systems of geometric paths. We start with critical
percolation, and then consider two other discrete models with
particular harmonic properties, the Loop-Erased Random Walk (LERW) and the
Uniform Spanning Tree (UST). We then proceed to show how to use the
Poissonian structure of the loop-soup to define a natural notion of
non-intersection for $n$ chordal $\SLE_\kappa$'s, where $0<\kappa\leq
8/3$.

\subsection{Crossing events for critical percolation}

We discuss here consequences of the locality property for
$\SLE_6$/critical percolation, in particular regarding holonomic
systems. 

\begin{figure}[htbp]
\begin{center}
\centerline{\psfig{file=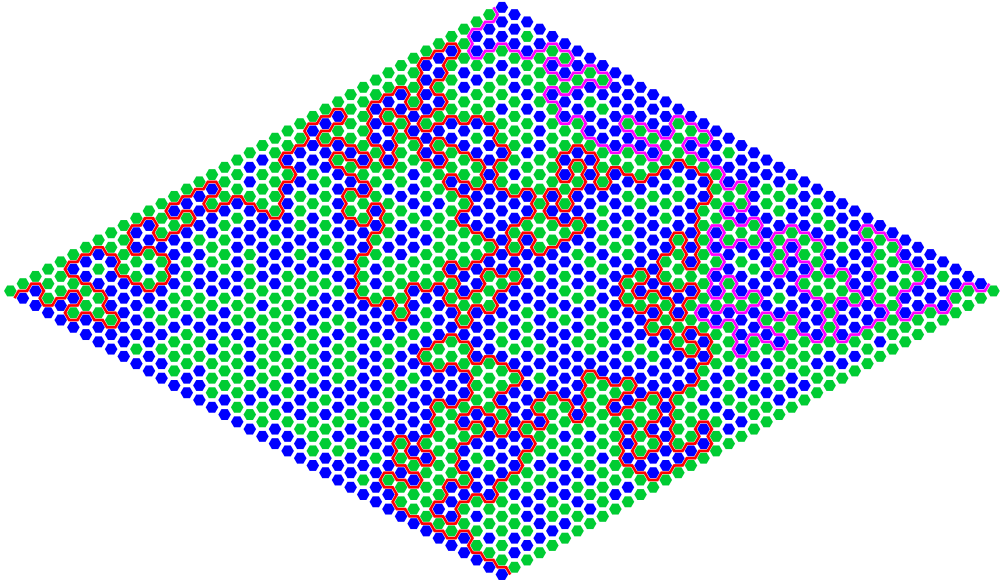}}
\end{center}
\caption{Two interfaces (red, fuchsia) bound a crossing of the lozenge}\label{figloz}
\end{figure}

First, we describe a family of percolation events that appear as
particularly well-suited for an $\SLE$ analysis. Recall the set-up of
Cardy's formula, relating to critical percolation in a conformal
quadrilateral. For simplicity, we will consider site percolation on
the triangular lattice: each site is colored in blue or yellow with
probability $1/2$, all sites being independent. (Alternatively, the
hexagons of a honeycomb tiling are colored in blue, yellow). A
subgraph of the
triangular lattice, with mesh $\eps\searrow 0$, approximates a
quadrilateral. As boundary conditions, a pair of opposite edges of the
quadrilateral is set to blue, while the other pair is set to
yellow. For (plane) topological reasons, either the two yellow edges
are connected by a yellow path on the lattice, or the two blue edges
are connected by a blue path (see Figure \ref{figloz}).

It seems quite natural to generalize
this construction to $(2n)$-gons, $n\geq 2$. So consider a bounded
Jordan domain $D\subset\C$; $a_1\dots a_{2n}$ are distinct points on
the boundary (in counterclockwise order, say). The boundary conditions
are alternate: for $1\leq i\leq n$, the arc $(a_{2i},a_{2i+1})$ is set to yellow, and
$(a_{2i-1},a_{2i})$ is set to blue (with cyclical indexing,
i.e. $a_{2n+1}=a_1$). We are interested in the connectivity properties
of the random percolation graph that are observable from the boundary.
Denote by $e_i$ the edge $(a_i,a_{i+1})$. Let $c(e_i,e_j)=1$ if
$e_i,e_j$ are of the same color and are connected by a path of this
color, and $c(e_i,e_j)=0$ otherwise. Note that a path can
include a portion of the boundary. An elementary event is
an event of type
$$C(\epsilon)=\bigcap_{1\leq i<j\leq 2n}\{c(e_i,e_j)=\epsilon_{i,j}\}$$
where $\epsilon=(\epsilon_{i,j})_{1\leq i<j\leq 2n}$,
$\epsilon_{i,j}\in\{0,1\}$. We now describe the non-empty elementary
events. We suppose that the mesh $\eps$ is small enough, so that the
edges are nonempty and disjoint. From the Russo-Seymour-Welsh theory
(and FKG inequality, see \cite{Grim}), the elementary events that are
topologically possible will happen with probability bounded away from
0 (and from 1) as the mesh $\eps$ goes to zero.

\begin{figure}[htbp]
\begin{center}
\input{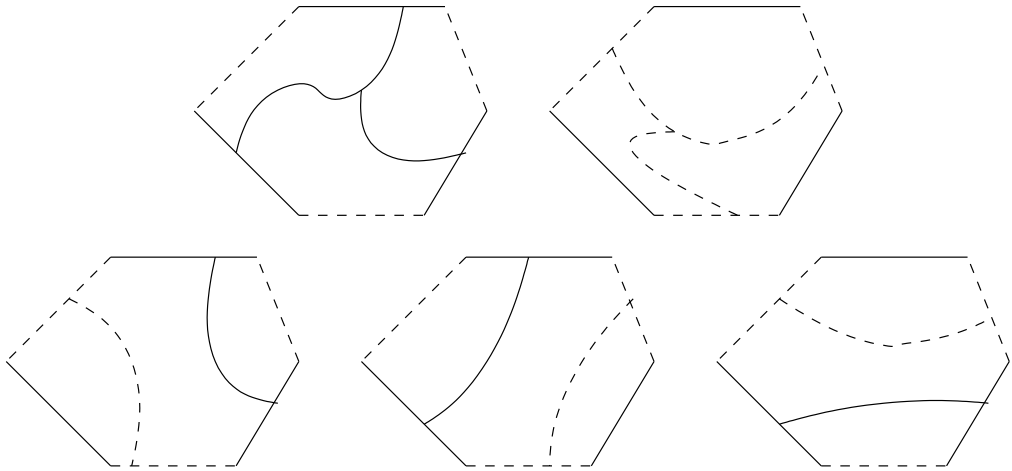}
\end{center}
\caption{The five configurations for a conformal hexagon (schematic)}\label{hexconfig}
\end{figure}

Let $n\geq 2$. There are exactly $C_n$ (non-empty) elementary events with
probability bounded away from 0 as $\eps\searrow 0$, where $C_n$ denotes the
$n$-th Catalan's number:
$$C_n=\frac{\binom{2n}{n}}{n+1}.$$
Indeed, it is easy to see that there is a one-to-one correspondence
between elementary events and non-crossing partitions of the set of
blue edges. As is well known, this number is $C_n$. This can be seen
as follows: consider the smallest $i>1$ such that $(a_{2i-1},a_{2i})$ is
connected to $(a_1,a_2)$ by a blue path; this leads directly to the
recurrence relation for Catalan's numbers. For instance, for $n=3$, one gets $C_3=5$ possible configurations for
a conformal hexagon with alternate boundary conditions, as illustrated
by Figure \ref{hexconfig}.

\begin{figure}[htbp]
\begin{center}
\scalebox{.25}{\input{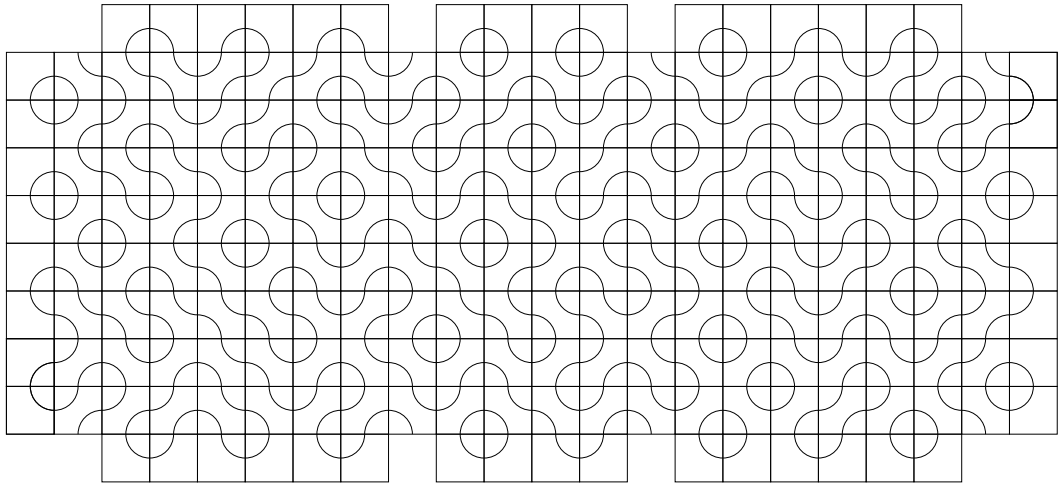}}
\end{center}
\caption{Four paths in bond percolation on $\Z^2$}\label{bondz2}
\end{figure}

Alternatively, at each point $a_i$ where boundary conditions change,
one can start an exploration path (\cite{S0,Sm1}) that winds between
the connected components of $(a_{i-1},a_i)$ and $(a_i,a_{i+1})$. Such a
path is simple and ends at $a_j$, for some $j$ ($i$ and $j$ have opposite parity). So the $(2n)$ boundary points are paired by $n$
non-intersecting exploration processes. This pairing is random; the
number of such pairings (satisfying the non-crossing condition) is
$C_n$. Note also that a site that touches two exploration
processes is pivotal for these events. For bond percolation on $\Z^2$,
it is convenient to represent the interfaces between connected
components in $\Z^2$ and connected components for the associated
percolation configuration on the dual graph. The collection of these
interfaces (with appropriate boundary conditions) is made of closed
loops and $n$ non-intersecting simple paths connecting the $(2n)$
``free'' points on the boundary (see Figure \ref{bondz2}).

The exploration process for critical site percolation on the
triangular lattice converges to
$\SLE_6$ in the scaling limit (\cite{Sm1,CamNew}). It is conjectured,
and supported by numerical evidence, that it is also the case for
critical bond percolation on $\Z^2$ (and also for more general lattices).

From conformal invariance of the scaling limit of critical
percolation, the probability of any of these elementary events should
define a function on the corresponding moduli space (i.e. the space of
Jordan domains with $(2n)$ (distinct) marked boundary points modulo conformal
equivalence). Denote by ${\mc M}_{2n}$ this moduli space, which can be
seen as a smooth $(2n-3)$-dimensional manifold, for $n\geq 2$ (for
a discussion of $\SLE$ and moduli space, see e.g. \cite{FK}). Considering all the
possible configurations, one defines a function ${\mc
  M}_{2n}\rightarrow \R^{C_n}$. Assuming that this function is smooth,
the dimension of the image is at most $(2n-3)$, which implies the
existence of a large number of smooth relations between probabilities
of elementary events. Indeed, the dimension of the moduli space,
$(2n-3)$, is negligible compared with $C_n$:
$$C_n=\frac{\binom{2n}{n}}{n+1}\sim
\frac 1n.\frac{\sqrt{2\pi 2n}(2n)^{2n}e^{-2n}}{(\sqrt{2\pi
    n}(n^n)e^{-n})^2}\sim \sqrt{\frac2\pi}.\frac {4^n}{n^{3/2}}$$
as follows from Stirling's formula. The sum of all probabilities is 1;
the nature of other smooth relations between these probabilities is
unclear. Letting a given edge of the conformal $(2n)$-gon shrink to 0,
one gets a conformal $(2n-2)$-gon, so that the corresponding part of
the boundary of ${\mc M}_{2n}$ can be identified with ${\mc
  M}_{2n-2}$. Considering this operation for any edge, it appears that
any affine relation between the $C_n$ probabilities is proportional to
the trivial (normalization) relation.

Assuming conformal invariance, one can give a differential
characterization of these probabilities. More specifically, let
$f(x_1,\dots x_{2n})$ be the probability (in the scaling limit) of
any one of the $C_n$ elementary events associated with the configuration
$(\H,x_1\dots x_{2n})$, $x_1<\cdots<x_{2n}$. Let $i\in\{1\dots 2n\}$;
we consider an infinitesimal percolation hull at $x_i$. The boundary changes
color at $x_i$; Smirnov's key result (see \cite{Sm1}) is that the
percolation exploration process started from $x_i$ converges to $\SLE_6$ started
at $x_i$ (be it chordal $\SLE_6$ to another boundary point, or radial
$\SLE_6$ to an inner point, since all these are equivalent for short
enough times). Let $(\gamma_u)_{0\leq u\leq\eta}$ be the exploration
process, and $(K_u)_{0\leq u\leq \eta}$ be the associated family
of hulls, where $\eta$ is such that the vertices $x_j$, $j\neq i$, are
not disconnected at time $\eta$. Then the elementary event holds for
$(\H,x_1,\dots, x_{2n})$ if and only if it holds for the (random) conformal
$(2n)$-gon $(\H\setminus K_u,x_1,\dots,\gamma_u,\dots x_{2n})$, for any
$u\in [0,\eta]$;
assuming conformal invariance, this has probability $f(g_u(x_1),\dots,
W_u,\dots g_u(x_{2n}))$, where $(g_u)$ are the conformal equivalences
defining the $\SLE_6$ process. So $(f(g_u(x_1),\dots,
W_u,\dots g_u(x_{2n})))_{0\leq u\leq\eta}$ is a martingale. To sum up,
assuming that the probability of an elementary event defines a smooth
function on ${\mc M}_{2n}$, then the function $f$ is annihilated by the
following differential ideal:
$$I_n=\langle {\mc L}_1,\dots {\mc L}_{2n},\ell_{-1},\ell_0,\ell_1\rangle$$
where ${\mc L}_i=3\partial_{ii}+\sum_{j\neq i}\frac
2{x_j-x_i}\partial_j$ (infinitesimal generator for $\SLE_6$ growing at
$x_i$), and $\ell_k=-\sum x_i^{1+k}\partial_i$. Note that the (real) Lie algebra
generated by the $(\ell_k)_{k\in\Z}$ is (isomorphic to the) Witt algebra, and
that the subalgebra $\langle \ell_{-1},\ell_0,\ell_1\rangle$ is
isomorphic to $\mathfrak {sl}_2(\R)$, the tangent algebra of the
Moebius group (the group of conformal automorphisms of a simply connected
domain). The fact that one can explore the $(2n)$-gon starting at any of its
vertices is a feature of locality.

It is readily seen that $I_n$, which is an ideal of differential
operators with rational coefficients, is holonomic, i.e. the vector space of
functions annihilated by $I_n$ in the neighbourhood of a generic point
has finite dimension (the rank of $I_n$ is the dimension of this
space). For background on holonomic systems, see e.g. \cite{Yosh}.
It is elementary that the rank of $I_n$ is no greater
that $4^n$. Indeed, if $f$, defined in a neighbourhood $U$ of a generic
point, is annihilated by $I_n$, let $F:U\rightarrow\R^{\{0,1\}^{2n}}$,
where:
$$F=(\partial_1^{\epsilon_1}\dots\partial_{2n}^{\epsilon_{2n}}f)_{\epsilon\in
  \{0,1\}^{2n}}.$$   
Using the operators ${\mc L}_i$, it readily follows that one can write
$\partial_i F=M_iF$ for $i=1\dots 2n$, where $M_i$ is a matrix of
rational coefficients. A local solution of this system in a neighbourhood of
$(x_1,\dots,x_{2n})$ is entirely determined by the initial condition
$F(x_1,\dots x_{2n})$, so the solution space is of dimension at most
$2^{2n}=4^n$ (it would be exactly of dimension $4^n$ if the Frobenius integrability
conditions were satisfied: $\partial_jM_i-\partial_iM_j=M_jM_i-M_iM_j$). Since $F_{(0,\dots, 0)}=f$, the rank of $I_n$ is at most
$4^n$. Note that we have not used the operators $\ell_{-1}$,
$\ell_0$, $\ell_1$, so this is a crude estimate. 
Alternatively, it is easily seen that the dimension of the
characteristic variety of $I_n$ is $(2n)$, implying holonomy (see
e.g. \cite{Yosh}). Making use of
conformal invariance, if $f$ is annihilated by $I_n$, then one can write:
$$f(x_1,\dots, x_{2n})=g([x_1,x_2,x_3,x_4],\dots,[x_1,x_2,x_3,x_{2n}])$$
where $g$ is a function of the $(2n-3)$ variables
$y_i=[x_1,x_2,x_3,x_{3+i}]$, where $[.,.,.,.]$ is the
cross-ratio. Since ${\mc L}_4,\dots {\mc L}_{2n}$ annihilate $f$, $g$ is
annihilated by operators $\tilde {\mc L}_1,\dots \tilde {\mc L}_{2n-3}$ such that
the only second order differential term featuring in $\tilde {\mc L}_i$ is
$\partial^2/\partial y_i^2$. Reasoning as above, this implies that the
rank of $I_n$ is at most $2^{2n-3}$. We expect the rank to be exactly $C_n$.
 
In the case $n=2$, if a function $f$ is annihilated by $I_2$, then by
conformal invariance ($\ell_{-1}f=\ell_0f=\ell_1f=0$), one can write
$$f(x_1,\dots,x_4)=g([x_1,\dots, x_4])$$
where $[.,.,.,.]$ is the cross-ratio. Then $f$ annihilates $I_2$ if
and only if $g$ satisfies the second order ODE:
$$3u(1-u)g''(u)+(2-4u)g'(u)=0$$
which is a Fuchsian differential equation with three singular regular
points on the Riemann sphere (namely $0$, $1$, $\infty$), so it is
essentially equivalent to a hypergeometric differential equation (see
\cite{Yosh}, p.26). 
It is easily seen that the rank of $I_2$ is exactly $2$; at a generic
point, the solution space is generated by a constant function and by
the function appearing in Cardy's formula. It is worth remarking that we
could have chosen any cross-ratio (there are six of them) to get the same
solution space, which translates into exceptional invariance properties
of the considered ODE (under homographies that fixate the singular
locus $\{0,1,\infty\}$). In fact, this is the only second-order Fuchsian ODE (with
singular locus $\{0,1,\infty\}$) that is invariant under the
substitutions $u\mapsto
(1-u)$, $u\mapsto 1/u$, so that the three singular points play exactly
symmetric roles.

\subsection{Loop-Erased Random Walks and Fomin's formulae}

In \cite{Fomin}, Fomin considers a loop-erased version of the well-known Karlin-McGregor Formula.
Problems pertaining to the scaling limits of Fomin's formulae are
studied in \cite{KozL}. From Wilson's algorithm (\cite{Wi}), we know that there is
an exact identity between branches of a uniform spanning tree (UST)
and loop-erased random walks (LERW). It turns out that the natural way
to define non-intersecting LERWs is to consider disjoint branches of
an ambient UST.

More precisely, consider the following situation. A simply connected
domain $D$ of the plane, with, say, smooth Jordan boundary, is
approximated by a subgraph of a lattice with mesh $\eps$. Then take a
uniform spanning tree of this graph with wired boundary
conditions. Consider now $x_1,\dots x_n$ $n$ points on the boundary
and $y_1,\dots,y_n$ $n$ points one lattice spacing away from the
boundary, such that $x_1,\dots,x_n,y_n,\dots y_1$ are in
counterclockwise order. We condition on the event that the minimal
subtree containing $y_1,\dots,y_n$ and the boundary has no triple
point in the bulk, and that the branch of $y_i$ connects to the
boundary at $x_i$, $i=1,\dots,n$. Then take the scaling limit of these
discrete paths. We get $n$ non-intersecting paths connecting $x_i$ to
$y_i$, $i=1,\dots,n$; each of them has density w.r.t. chordal $\SLE_2$
(\cite{LSW2}).

In the discrete setting, let $H$ be the harmonic measure for simple
random walk on the lattice killed when it reaches the boundary. Then
the probability of the event considered above is given by Fomin's
formula:
$$\det(H(y_i,\{x_j\}))_{1\leq i,j\leq n}.$$
As follows from \cite{KozL}, in the scaling limit (say in the upper
half-plane $\H$), the Girsanov density of the system of paths w.r.t.
independent chordal $\SLE_2(x_i\rightarrow y_i)$ is given by the
scaling limit of this determinant divided by the product of its
diagonal terms:
$$\det\left(\frac 1{(x_i-y_j)^2}\right)_{i,j}\prod_i (x_i-y_i)^2.$$ 
In the continuum limit, we will see later an extension of this
situation to values of $\kappa$ between $0$ and $8/3$. Note that this
discrete construction is reversible
(for slightly different boundary conditions). This situation 
further illustrates the close relationship between LERW and discrete
harmonic measure and its application to scaling limits. In the next
section, we shall make use of the connection  between UST and
reflected harmonic measure.

\subsection{Uniform Spanning Trees}

If we specialize the results we shall obtain later in this paper to $\kappa=8$, the formulae
become symmetric in the $(2n)$ variables and have a nice Riemann
surface interpretation. On the other hand, chordal $\SLE_8$ is known
to be the scaling limit of (the Peano curve of) the Uniform Spanning
Tree (see \cite{LSW2}). It is not hard to think of boundary conditions
that enable to define multiple exploration processes in the discrete
setting. From Wilson's algorithm (\cite{Wi}), it is then possible to extract the
drift terms. Though, the connection with our formulae is not quite
immediate, so we shall give some details here.

\begin{figure}[htbp]
\begin{center}
\scalebox{.8}{\input{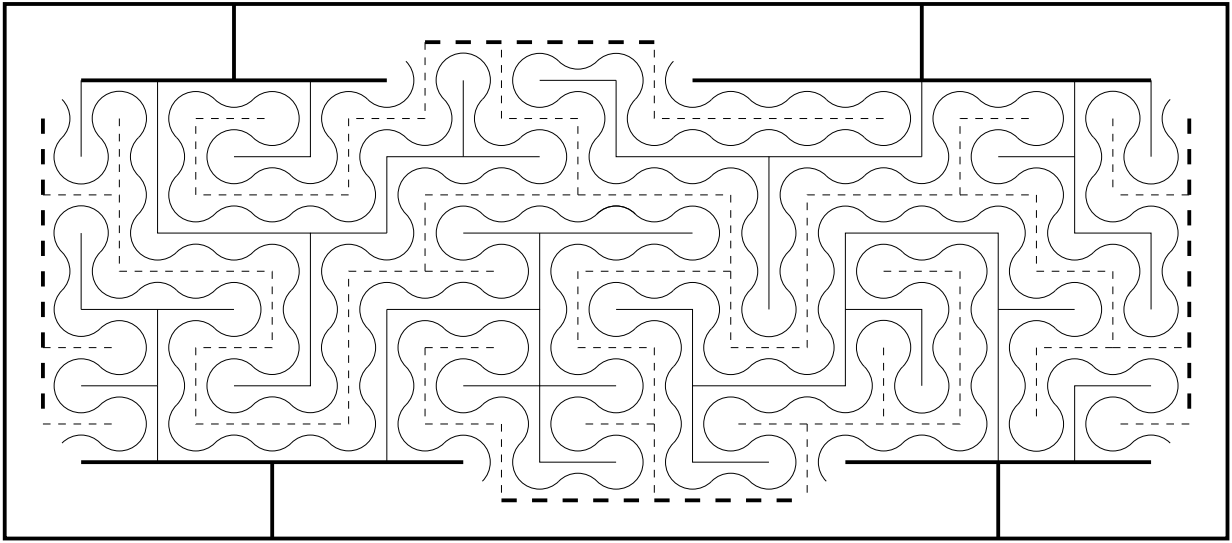}}
\end{center}
\caption{Tree (solid), dual tree (dashed), Peano paths (curved)}\label{ust2}
\end{figure}

Consider a discrete lattice approximation of a bounded simply
connected with, say, piecewise smooth boundary and $(2n)$ marked
points $x_1,\dots x_{2n}$ on the boundary. The $n$ boundary arcs
$(x_2,x_3),(x_4,x_5),\dots, (x_{2n},x_1)$ are wired, and connected
by an external wiring; the $n$ other boundary arcs are free (i.e. wired
for the dual graph). One then samples uniformly from spanning trees
with these conditions. If we erase the external wiring, one gets
$n$ trees (corresponding to the $n$ wired boundary components). Each vertex
$z$ is connected to one of the wired boundary component with
probability given by the (graph) harmonic measure, from Wilson's
algorithm. In the scaling
limit, this should converge to harmonic measure with normal reflection
on the free parts (see Figure \ref{ust2}). 

With these conditions, the Peano path starting at $x_{2i}$ ends
at $x_{2i+1}$, so that the boundary conditions fix a pairing of the $(2n)$
points; note also that this pairing is different from the one most
natural for Fomin's formulae. Assume now for simplicity that the
domain in the upper half-plane $\H$, and $x_1<x_2<\cdots<x_{2n}$.
We want to prove that this situation
corresponds in the scaling limit to $\SLE_8$'s whose drift terms are
given by log-derivatives of the ``partition function'':
$$\psi({\bf x})=\prod_{1\leq i<j\leq 2n}(x_j-x_i)^{1/4}\int_C \prod_{1\leq i<j<n}(u_j-u_i)\prod_{i=1}^n\frac {du_i}{\left(\prod_{j=1}^{2n}(u_i-x_j)\right)^{1/2}}$$
where $C$ is a Cartesian product of $C_1,\dots, C_{n-1}$, where $C_i$ is
a cycle circling clockwise round the
segment $(x_{2i-1},x_{2i})$. Note that one can choose a determination
of the integrand on such cycles. In this expression, rewrite the product
$\prod_{i<j}(u_j-u_i)$ as a Vandermonde determinant. Define also:
$$\omega_i=\frac{u^{i-1}du}{\left(\prod_{j=1}^{2n}(u-x_j)\right)^{1/2}}, \omega'_i=\frac{(u-x_1)^{i-1}du}{\left(\prod_{j=1}^{2n}(u-x_j)\right)^{1/2}}$$
for $i\in\{1,\dots n-1\}$. It is well-known that
$(\omega_1,\dots,\omega_{n-1})$ is a basis of abelian differentials of
the first kind for the hyperelliptic curve of genus $g=n-1$:
$$t^2=(s-x_1)\cdots (s-x_{2n})$$
and $(C_1,\dots,C_{n-1})$ is half of a canonical homology basis for
this curve (see e.g. \cite{FarKra} for background on compact Riemann
surfaces). Then $\psi$ can be written as a determinant of periods:
$$\psi({\bf x})=\prod_{1\leq i<j\leq 2n}(x_j-x_i)^{1/4}\det\left(\int_{C_j}\omega_i\right)=\prod_{1\leq i<j\leq 2n}(x_j-x_i)^{1/4}\det\left(\int_{C_j}\omega'_i\right).$$
Consider the Peano exploration starting from $x_1$ and ending at
$x_{2n}$. We want to recover the form of $\psi$ from discrete
arguments. First, we note that the situation can be seen as a
degenerate case of a UST in a multiply-connected domain. Consider
$n-1$ compact holes $K_1,\dots K_{n-1}$ in $\H$, with the following conditions: two points
$x_1$ and $x_{2n}$ are marked on $\R$; $(x_1x_{2n})$ is free and
$(x_{2n}x_1)$ is wired; the boundaries of the holes are wired, and all
the wired parts are considered as wired together. One recovers the
previous situation is the holes are segments close to $\R$. If $z$ is
any bulk point, then Wilson's algorithm provides $n-1$ martingales for
the Peano exploration process. More precisely, let $H_t$ be the
remaining domain after time $t$ of the exploration (for some time
parameterization); and ${\rm Harm}$ is the harmonic measure in $H_t$
with normal reflection on the free part of the boundary
($(x_1x_{2n})$). Then:
$$t\mapsto M_t=\vphantom{H}^t({\rm Harm}(z,\partial K_1),\dots {\rm Harm}(z,\partial K_{n-1}))$$
is a vector-valued martingale (see \cite{LSW2}). Following
\cite{LSW2}, and the work of Makarov and
Zhan for SLE in multiply connected domains (see \cite{DZT}), one determines from this
the driving process of the scaling limit. We use here the notations
and conventions of \cite{D6}. So $(g_t)$ is a family of conformal
equivalences $H_t\rightarrow \H$ extending through the holes and with
hydrodynamic normalization at infinity:
$$\partial_t g_t=\frac {2}{g_t-X_t}, dX_t=\sqrt \kappa
dB_t+\kappa\frac {\partial_x\psi}{\psi}dt$$
for some covariant function $\psi$ of the configuration, which we want
to identify. Let $z$ be a point on the free part of the boundary
close to $x=x_1$; expanding at $z=x$, one gets:
$$M_t=(z-x)^{1/2}a_t+(z-x)^{3/2}b_t+\cdots$$
the half-integer exponents coming from the normal reflection. The martingale
condition at $t=0$ gives:
$$\left(\frac \kappa 2\partial_{xx}+\kappa\frac {\partial_x\psi}\psi\partial_x+\frac{2}{z-x}\partial_z+\cdots\right)\left((z-x)^{1/2}a+(z-x)^{3/2}b+\cdots\right)=0$$
where $a$ and $b$ are now (vector-valued) functions of the initial
configuration. Considering the terms in $(z-x)^{-3/2}$ and
$(z-x)^{-1/2}$ in this equation, one gets the necessary conditions
$\kappa=8$ and:
$$4\frac {\partial_x\psi}{\psi} a+4\partial_x a=6b.$$
So the drift term $\partial_x\log\psi$ can be recovered from $a$, $b$, that
are coefficients in the expansion of the harmonic measure. Note that
$\partial_x\log\psi$ is a real number, and the condition is a relation
on vectors. Now we want to
prove that $\psi$ as defined above satisfies this equation in the
degenerate case where $K_1=(x_2x_3),\dots
K_{n-1}=(x_{2n-2}x_{2n-1})$. We assume that $n>2$ (the case $n=2$
being comparatively trivial).
 
Let $P$ be the matrix:
$$P=\left(\int_{C_i}\omega'_j\right)_{1\leq i,j<n}$$
and $Q(z)$ be the vector
$\int_x^z\vphantom{()}^t(\omega'_1,\dots,\omega'_{n-1})$. Then it is easy to see
that $\Re(P^{-1}Q)$ is the image of the harmonic measure vector $H_0$
by a fixed triangular matrix. In fact, the imaginary part will also give a
martingale. So we can replace the harmonic measure vector by $P^{-1}Q$ in
what follows, for simplicity. Hence, if $(e_1,\dots e_{n-1})$ is the standard basis of
$\C^{n-1}$, one gets:
$$a=P^{-1}(2e_1)\prod_{i>1}(x_i-x)^{-1/2},
b=P^{-1}\left(-\frac 13\left(\sum_{i>1}\frac
1{x-x_i}\right)e_1+\frac 23 e_2\right)\prod_{i>1}(x_i-x)^{-1/2}$$
where $P^{-1}Q(z)=(z-x)^{1/2}a+(z-x)^{3/2}b+\cdots$. Using the fact
that:
$$\partial_x\det(P)=\det(P)\Tr(P^{-1}\partial_xP)$$
and $\psi=\prod_{1\leq i<j}(x_j-x_i)^{1/4}\det(P)$, after
simplifications, we have just to
check that:
$$8\Tr(P^{-1}\partial_xP)e_1+8P\partial_x P^{-1}e_1=4e_2.$$
This follows from $P\partial_x (P^{-1})+(\partial_x P) P^{-1}=0$ and the
fact that $\partial_x \omega'_{i+1}=(1/2-i)\omega'_i$ if
$i\in\{1,\dots,n-2\}$, so that:
$$\partial_x P=\left(
\begin{array}{ccccc}
\ast&\ast&\cdots&\cdots&\ast\\
-\frac 12&0&\dots&\dots&0\\
0&\ddots&\ddots\\
\vdots&\ddots&\ddots&\ddots\\
0&\dots&0&\frac 52-n&0
\end{array}\right)P$$
and the previous identity can be read from the first column of $(\partial_x P) P^{-1}$.

\subsection{Non-intersection and the restriction property}

As pointed out in \cite{D6}, the restriction property for $\SLE_{8/3}$
can be put to use to get natural constructions of non-intersecting
SLEs. For simplicity, consider first a simply-connected domain in $\C$
with four points marked on the boundary, say $(\H,x_1,x_2,x_3,x_4)$,
where $\H$ is the upper half-plane. Consider two independent
$\SLE_{8/3}$'s in $\H$, from $x_1$ to $x_2$ and $x_3$ to $x_4$ resp.,
with traces $\gamma,\gamma'$. Say $(g_t)$ is the family of conformal
equivalences associated with the first one (time parameterization is
unimportant here). 
The restriction property for the second $\SLE$
implies that the law of $g_t(\gamma')$ conditionally on
$\{\gamma\cap\gamma'=\varnothing\}$ is the law of an $\SLE_{8/3}$ from
$g_t(x_3)$ to $g_t(x_4)$ conditioned not to intersect
$g_t(\gamma_{(t,\infty)})$. Using at the same time the Markov property
and the restriction property for the first $\SLE$, one gets that
$(g_t(\gamma),g_t(\gamma'))$ conditionally on
$\{\gamma\cap\gamma'=\varnothing\}$ is distributed as
$(\tilde\gamma,\tilde\gamma')$, where $\tilde\gamma$ and $\tilde\gamma'$ are independent $\SLE_{8/3}$'s
going from $g_t(\gamma_t)$ to $g_t(x_2)$ and from $g_t(x_3)$ to $g_t(x_4)$
respectively, conditioned not to intersect. Indeed:
\begin{align*}
{\mc L}(g_t(\gamma_{(t,\infty)}),g_t(\gamma')|\gamma\cap\gamma'=\varnothing)
&={\mc
L}(g_t(\gamma_{(t,\infty)}),g_t(\gamma')|\gamma_{(0,t)}\cap\gamma'=\varnothing,g_t(\gamma_{(t,\infty)})\cap
g_t(\gamma')=\varnothing
)\\
&={\mc
L}(\tilde\gamma,\tilde\gamma'|\tilde\gamma\cap
\tilde\gamma'=\varnothing)
\end{align*}
Moreover, the system $(\gamma,\gamma')$ has also a restriction
property. 
This can be generalized to $n$ SLEs connecting $(2n)$ points on the
boundary, and conditioned on no intersection between
any two of them. This induces a topological constraint on the order of
the marked points $(x_1,\dots,x_{2n})$ on the boundary. If the points
$(x_1,\dots,x_{2n})$ are in cyclical order on the boundary, one is
interested in ``non-crossing pairings'' of these points (deciding that
an SLE connects $x_i$ and $x_j$ pairs $x_i$ and $x_j$). The number of
such pairings is given by Catalan's number $C_n$. The value of the
probability of non-intersection associated with such a pairing can be
seen as a function on the moduli space ${\mc M}_{2n}$, and satisfies
$(2n)$ evolution equations (\cite{D6}).

If $\kappa\in (0,8/3)$, this can be extended
using loop-soups of intensity
$\lambda_\kappa=(6-\kappa)(8-3\kappa)/2\kappa$ (see \cite{LSW3,LW}). More
precisely, consider $(2n)$ marked points on the boundary, $n$
independent $\SLE_\kappa$ connecting these points, and independent
loop-soups $L_2,\dots,L_n$ with intensity $\lambda_\kappa$. If the pairing is
non-crossing, there is a positive probability that no loop in $L_j$
intersects $\cup_{i<j}\gamma^i$ and $\gamma^j$, $j=2,\dots,n$.
Conditioning on this event (and
assuming reversibility), one gets a Markov property as above for
each of the boundary points. Here one uses the Markov property of SLE,
the restriction property for the loop-soup, the identity connecting
loop-soup, SLE, and restriction measures (\cite{LSW3}), and 
the Poissonian nature of the loop-soup. Note that one can rephrase the
conditioning in a symmetric fashion. Indeed, if $\nu$ is the Brownian
loop measure (denoted by $\mu^{loop}$ in \cite{LW}), conditionally on
the position of the SLEs, the probability that no loop in $L_j$
intersects $\cup_{i<j}\gamma^i$ and $\gamma^j$, $j=2,\dots,n$ is:
$$\exp\left(-\lambda_\kappa\sum_{j}\nu(\{\delta: \delta\cap(\cup_{i<j}\gamma^i)\neq\varnothing,\delta\cap\gamma^j\neq\varnothing\})\right).$$
It readily appears that the sum in the exponent can be written as:
$$\sum_{n>0}\sum_{i_1<\cdots< i_n}(n-1)\nu(\{\delta:
\delta\cap\gamma^j\neq\varnothing {\rm\ iff\ }
j\in\{i_1,\dots,i_n\}\})$$
by an inclusion-exclusion argument.
So it is equivalent to condition on no loop in $L_j$ intersecting (at
least) $j$
distinct SLEs, $j=2,\dots,n$.

Let us now formulate precisely and prove these results.
Let $\kappa\in (0,8/3]$. In the domain $(\H,x_1,y_1,\dots,x_n,y_n)$, consider $n$ independent
$\SLE_\kappa$'s from $x_i$ to $y_i$, $i=1,\dots,n$, and independent
loop-soups $L_2,\dots,L_n$ with intensity $\lambda_\kappa$. The SLEs are defined by
conformal equivalences $(g^i_t)$ and traces $\gamma^i$.
Assuming that the
pairing $(x_1,y_1),\dots,(x_n,y_n)$ is non-crossing, denote by $\mu_{(x_1,\dots, y_n)}$ the
law of the SLEs conditionally on the event that no two SLEs intersect
and no loop in $L_j$ intersects $\gamma^j$ and $\cup_{i<j}\gamma^i$ (if $\kappa<8/3$, one can drop
the first condition). A hull $A$ is a compact subset
of $\overline\H$ such that $\H\setminus A$ is simply connected,
$A\cap\R\subset\overline{A\cap\H}$ and
$x_1,\dots, y_n$ are not in $A$; $\phi_A$ is a conformal equivalence
$\H\setminus A\rightarrow A$. Let $L$ be yet another independent loop-soup
with intensity $\lambda_\kappa$; the union of $A$ and loops that
intersect it is denoted by $A^{L}$; $L^A$ is the collection of loops in
$L$ that do not intersect $A$.

\begin{Prop}\label{Prestr}
Under the previous assumptions:\\
(i) (Markov property) Under  $\mu_{(x_1,\dots, y_n)}$, the law of
$g^i_t(\gamma^1,\dots,\gamma^i_{[t,\infty)},\dots \gamma^n)$ is that
of $(\gamma^1,\dots,\gamma^n)$ under $\mu_{(g^i_t(x_1),\dots g^i_t
(\gamma^i_t),\dots, g^i_t(y_n))}$ (up to time reparameterization).\\
(ii) Let $\psi(x_1,\dots,y_n)$ be the probability that no two SLEs intersect
and no loop in $L_j$ intersects $\gamma^j$ and $(\cup_{i<j}\gamma^i)$,
$j=2,\dots,n$. Let $g_t=g^i_t$ for some
$i\in\{1,\dots,n\}$. Then the process:
$$\psi(g_t(x_1),\dots,g_t(\gamma^i_t),\dots
g_t(y_n))\prod_{j\neq i}\left(g'_t(x_j)g'_t(y_j)\left(\frac{y_j-x_j}{g_t(y_j)-g_t(x_j)}\right)^2\right)^{\alpha_\kappa}$$  
is a martingale under the (unconditional chordal) SLE measure for the $i$-th SLE
, where $\alpha_\kappa=(6-\kappa)/2\kappa$.\\
(iii) (Restriction property) Under $ \mu_{(x_1,\dots,
y_n)}$, the probability that $A^{L}$ does not intersect
$(\cup\gamma^i)$ is given by:
$$\frac{\psi(\phi(x_1),\dots,\phi(y_n))}{\psi(x_1,\dots,y_n)}\prod_i \left(\phi'(x_i)\phi'(y_i)\left(\frac{x_i-y_i}{\phi(x_i)-\phi(y_i)}\right)^2\right)^{\alpha_\kappa}$$
where $\phi=\phi_A$. Moreover, the image under $\phi$ of
$(\gamma_1,\dots,\gamma_n,L_2^A,\dots L_n^A,L^A)$ conditionally on 
``no loop in $L_j$
connects $\gamma_j$ to $\cup_{i<j}\gamma_i$, $j=2,\dots,n$, and no SLE is connected to $A$ by a loop in
$L$'' is distributed as $(\tilde\gamma_1,\dots,\tilde\gamma_n,\tilde
L_2,\dots, \tilde L)$ where $\tilde\gamma^i$ are independent SLEs
connecting $\phi(x_i)$ to $\phi(x_j)$, $\tilde L_2,\dots,\tilde L$ are
independent loop-soups, conditionally on ``no loop in $\tilde L_j$
connects $\tilde\gamma_j$ to $\cup_{i<j}\tilde\gamma_i$, $j=2,\dots,n$''.
\end{Prop} 

\begin{proof}
If $n=1$, (i) is the usual Markov property for SLE, (iii) is a result
from \cite{LSW3}, and (ii) is an empty statement (since
$\psi(x_1,y_1)=1$). 

By induction, assume that the results are proved for $(n-1)$
SLEs. 
Let us begin with (iii). Consider the event ``no loop in $L_j$
connects $\gamma_j$ to $\cup_{i<j}\gamma_i$, $j=2,\dots,n$, and no SLE is connected to $A$ by a loop in
$L$''. The probability of this event, conditionally on the $n$
independent SLEs is:
$$\exp\left(-\lambda_\kappa\left(\sum_{j\leq n}\nu(\{\delta: \gamma^j\cap\delta\neq\varnothing, (\cup_{i<j}\gamma^i)\cap\delta\neq\varnothing\})+\nu(\{\delta: (\cup_j\gamma^j)\cap\delta\neq\varnothing,A\cap\delta\neq\varnothing\})
\right)\right)
$$
For illustration, let us momentarily assume that $n=2$. Then the sum
in the exponential density reduces to:
$$\nu(\{\delta:\delta\cap\gamma^1\neq\varnothing,\delta\cap\gamma^2\neq\varnothing\})+\nu(\{\delta:\delta\cap(\gamma^1\cup\gamma^2)\neq\varnothing,\delta\cap
A\neq\varnothing\})$$
which we can rewrite by inclusion-exclusion as:
$$\nu(\{\delta:\delta\cap\gamma^1\neq\varnothing,\delta\cap
A\neq\varnothing\})+\nu(\{\delta:\delta\cap\gamma^2\neq\varnothing,\delta\cap
A\neq\varnothing\})+\nu(\{\delta:\delta\cap\gamma^1\neq\varnothing,
\delta\cap\gamma^2\neq\varnothing,
\delta\cap
A=\varnothing\}).
$$
The first term corresponds to the density of chordal $\SLE$ in
$(\H\setminus A,x_1,y_1)$ w.r.t. chordal $\SLE$ in $(\H,x_1,y_1)$;
symmetrically, the second term is the density of chordal $\SLE$ in
$(\H\setminus A,x_2,y_2)$ w.r.t. chordal $\SLE$ in
$(\H,x_2,y_2)$. This follows from the restriction property for a
single SLE, applied twice. The remaining term is the mass of loops in
$\H\setminus A$ that do not intersect the two SLEs. So (iii) follows
when $n=2$. Let us get back to the general case: $n\geq 2$, and we
assume that the assertions hold for $(n-1)$ SLEs.

The sum in the exponential can be written as:
\begin{align*}
\sum_{j<n}\nu(\{\delta: \gamma^j\cap\delta\neq\varnothing,
(\cup_{i<j}\gamma^i)\cap\delta\neq\varnothing\})\\
+\nu(\{\delta: \gamma^n\cap\delta\neq\varnothing,
(\cup_{i<n}\gamma^i)\cap\delta\neq\varnothing,A\cap\delta=\varnothing\})\\
+2\nu(\{\delta: \gamma^n\cap\delta\neq\varnothing,
(\cup_{i<n}\gamma^i)\cap\delta\neq\varnothing,A\cap\delta\neq\varnothing\})\\
+\nu(\{\delta:
\gamma^n\cap\delta\neq\varnothing,
(\cup_{i<n}\gamma^i)\cap\delta=\varnothing,A\cap\delta\neq\varnothing\})\\
+\nu(\{\delta:
\gamma^n\cap\delta=\varnothing,
(\cup_{i<n}\gamma^i)\cap\delta\neq\varnothing,A\cap\delta\neq\varnothing\})\\
\end{align*}
and after rearranging:
\begin{align*}
\sum_{j<n}\nu(\{\delta: \gamma^j\cap\delta\neq\varnothing,
(\cup_{i<j}\gamma^i)\cap\delta\neq\varnothing\})+\nu(\{\delta:
(\cup_{i<n}\gamma^i)\cap\delta\neq\varnothing,A\cap\delta\neq\varnothing\})\\
+\nu(\{\delta:
\gamma^n\cap\delta\neq\varnothing,A\cap\delta\neq\varnothing\})\\
+\nu(\{\delta: \gamma^n\cap\delta\neq\varnothing,
(\cup_{i<n}\gamma^i)\cap\delta\neq\varnothing,A\cap\delta=\varnothing\})
\end{align*}
The first line is treated by the induction hypothesis: 
it corresponds to the (unnormalized) density of $\mu_{(\H\setminus
A,x_1,\dots,y_{n-1})}$ w.r.t. $\mu_{(\H,x_1,\dots,y_{n-1})}$
The second line corresponds to the density of $\gamma^n$ as an SLE in
$(\H\setminus A,x_n,y_n)$ w.r.t. $\gamma^n$ as an SLE in the full
domain $(\H,x_n,y_n)$ (restriction property for a single SLE). The
last line is the mass of loops in $\H\setminus A$ that intersect both
$\gamma_n$ and $\cup_{i<n}\gamma^i$ (from the restriction property of
the loop-soup). So the exponential above is indeed the density of 
$\mu_{(\H\setminus
A,x_1,\dots,y_n)}$ w.r.t. $\mu_{(\H,x_1,\dots,y_n)}$. The statement on
the joint law with loop-soups follows from the restriction property
for the loop-soups $L_2,\dots,L$. (Note that we have exchanged loops
between loop-soups, but these were only loops intersecting $A$).

The formula in (iii) is obtained by keeping track of masses of
(unnormalized) measures in the induction scheme above; the new
covariance factor in the product comes from the restriction property
for a single SLE (say $\gamma^n$). 

If $A,B$
are two hulls and $A.B$ is the hull satisfying
$\phi_{A.B}=\phi_B\circ\phi_A$, then it is clear from the formula that
the probability that $\cup_i\phi_A(\gamma^i)$
does not touch $B^L$, conditionally on ``no loop in $L_j$
connects $\gamma_j$ to $\cup_{i<j}\gamma_i$, $j=2,\dots,n$, and no SLE is connected to $A$ by a loop in
$L$'' multiplied by the probability that $\cup_i\phi_A(\gamma^i)$
does not touch $A^L$
is the probability that $\cup_i\gamma^i$ does not
touch $(A.B)^L$. This gives the (usual) restriction property for
$(\cup_i\gamma^i)^L$, which is less precise than the statement (iii). 
(See also \cite{D6} for a discussion of these restriction measures).

Let us now prove (i),(ii). As explained above, the conditioning
is in fact symmetric in the $n$ SLEs. So, without loss of generality,
we can assume that $i=n$. Reasoning as above, using the Markov property for $\gamma^n$ and
the restriction property (iii) for $(\gamma^1,\dots,\gamma^{n-1})$
(where $A$ is replaced with $\gamma^n_{[0,t]}$),
one gets easily (i) and (ii).
\end{proof}

When $\kappa=2$, this situation connects with Fomin's determinantal
formulae (\cite{Fomin,KozL}). Indeed, adding the loops of a loop-soup
with intensity $\lambda_2$ to an $\SLE_2$, one gets a restriction
measure with exponent 1, that has the same outer boundary as the
Brownian Excursion. (In fact, the union of $\SLE_2$ with the unfilled loops that
intersect it has the same distribution as the range of the Brownian
excursion, see \cite{DZT,Law2}). For instance, consider the situation with two
$\SLE_2$'s. Then the probability that no loop intersects the two SLEs is
the same as the probability that the first SLE does not intersect the
second one with loops attached, which is the same as the probability
that an $\SLE_2$ does not intersect a Brownian excursion (with
respective endpoints $x_1,x_2,x_3,x_4$). This gives an interpretation
of the symmetry of Fomin's formulae in the scaling limit (this symmetry follows
from Wilson's algorithm in the discrete case).

Under mild regularity assumptions, It\^o's formula and (ii) imply that $\psi$, which is
conformally invariant by construction, is annihilated by the $n$ operators:
$$\frac\kappa 2\partial^2_{x_j}+\sum_{k\neq
j}\frac{2\partial_{x_k}}{x_k-x_j}
+\sum_{k}\frac{2\partial_{y_k}}{y_k-x_j}+\frac{2\partial_{y_j}}{y_j-x_j}+\frac{\kappa-6}{\kappa}\sum_{k\neq
j}\left(\frac{1}{x_k-x_j}-\frac{1}{y_k-x_j}\right)^2. $$
Under the assumption of reversibility (see \cite{RS01}), $\psi$ is also
annihilated by the $n$ operators obtained by swapping the $x$ and $y$
variables. Note that the pairing is materialized in the constant terms
of these operators. In fact, one can symmetrize the equations by
a conjugation. Indeed, denoting $x_{n+j}=y_j$, the
function
$$\psi(x_1,\dots,x_{2n})\prod_j(x_{n+j}-x_j)^{1-6/\kappa}$$
is annihilated by the operators
$$\left\{\begin{array}{l}\displaystyle\frac\kappa 2\partial_{kk}+\sum_{l\neq
k}\frac{2\partial_{l}}{x_l-x_k}+\frac{\kappa-6}\kappa\sum_{l\neq
k}\frac 1{(x_l-x_k)^2}, \hspace{1cm} k=1,\dots, 2n\\
\sum_k\partial_k\\
\sum_k x_k\partial_k-n(1-6/\kappa)\\
\sum_k x_k^2\partial_k-(1-6/\kappa)(x_1+\cdots+ x_{2n})
\end{array}\right.$$
the three first-order operators representing conformal covariance. 
This is the system (\ref{sys}).
So
under regularity and reversibility assumptions, each non-crossing
pairing of the $(2n)$ points $(x_1,\dots,x_{2n})$ gives rise to a
solution of this system.

In the case where $8/3<\kappa<4$, one can expect that one can still define a natural law on non-intersecting SLEs, whose density w.r.t. independent chordal SLEs is given by $\exp(-\lambda_\kappa\nu(\dots))$; this density is now unbounded:
\begin{align*}
d\mu_{(x_1,\dots,y_n)}(\gamma_1,\dots,\gamma_n)=&\psi(x_1,\dots,y_n)^{-1}\ind_{\{\gamma_i\cap\gamma_j=\varnothing,i<j\}}\exp\left(-\lambda_\kappa\sum_{j}\nu(\{\delta: \delta\cap(\cup_{i<j}\gamma^i)\neq\varnothing,\delta\cap\gamma^j\neq\varnothing\})\right)\\
& d\mu_{(x_1,y_1)}(\gamma_1)\dots d\mu_{(x_n,y_n)}(\gamma_n)
\end{align*}
where $d\mu_{(x_i,y_i)}$ is the measure induced on paths (say as
elements of the Hausdorff space) by standard chordal $\SLE_\kappa$. So
$\psi$ can be interpreted here as a (hopefully finite) 
partition function.

\section{Some particular cases} 

Let $\kappa>0$ and $n\geq 1$; $x_1,\dots, x_{2n}$ are boundary points
of $\H$. Consider the operators (\ref{sys}):
\begin{equation*}
\left\{\begin{array}{l}\displaystyle\frac\kappa 2\partial_{kk}+\sum_{l\neq
k}\frac{2\partial_{l}}{x_l-x_k}+\frac{\kappa-6}\kappa\sum_{l\neq
k}\frac 1{(x_l-x_k)^2}, \hspace{1cm} k=1,\dots, 2n\\
\sum_k\partial_k\\
\sum_k x_k\partial_k-n(1-6/\kappa)\\
\sum_k x_k^2\partial_k-(1-6/\kappa)(x_1+\cdots+ x_{2n})
\end{array}\right.
\end{equation*}
If $\varphi$ is annihilated by these operators and
$\iota$ is an involution of $\{1,\dots,2n\}$ with no fixed point (so
that $\iota$ determines a pairing), denote:
$$\psi(x_1,\dots,x_{2n})=\varphi(x_1,\dots,x_{2n})\prod_{\{j,\iota(j)\}}(x_j-x_{\iota(j)})^{6/\kappa-1}.$$
Then $\psi$ is conformally invariant (from the last three equations of
(\ref{sys})); moreover, consider an
$\SLE_\kappa$ from $x_j$ to $x_{\iota(j)}$, with associated conformal
equivalences $(g_t)$ (with hydrodynamic normalization) and driving
process $W_t$. Then:
$$\psi(g_t(x_1),\dots,W_t,\dots,g_t(x_{2n}))\prod_{\{k,\iota(k)\}\neq\{j,\iota(j)\}}\left(g'_t(x_k)g'_t(x_{\iota(k)})\left(\frac{x_k-x_{\iota(k)}}{g_t(x_k)-g_t(x_{\iota(k)})}\right)^2\right)^{\alpha_\kappa}$$
is a local martingale (from the $j$-th equation of (\ref{sys})).

We study this system in a few cases, before giving a set of formal solutions.

\subsection{Case $n=2$}
If $\varphi$ and $\psi$ are as above (say for the pairing $\{(x_1,x_2),(x_3,x_4)\}$), then $\psi$ is a conformally invariant function of
four boundary points, so  $\psi(x_1,\dots,x_4)=f(r)$, where $r$ is the cross-ratio:
$$r=\frac{(x_3-x_2)(x_4-x_1)}{(x_3-x_1)(x_4-x_2)}.$$
Now it is easy to check that $\varphi$ is a solution of (\ref{sys}) iff
$f(r)r^{-2/\kappa}$ is a solution of the hypergeometric equation with
parameters $(a,b,c)=(4/\kappa,1-4/\kappa,8/\kappa)$ (see \cite{B}).

Assume that $x_1,\dots,x_4$ are in cyclical order. It is easy to see
that the boundary conditions $\psi(x_1,x_3,x_3,x_4)=0$ ($r=0$),
$\psi(x_1,x_1,x_3,x_4)=1$ ($r=1$) determine the following solution:
$$\psi(x_1,x_2,x_3,x_4)=\frac{\Gamma(4/\kappa)\Gamma(12/\kappa-1)}{\Gamma(8/\kappa)\Gamma(8/\kappa-1)}r^{2/\kappa}\vphantom{F}_2F_1\left(\frac 4\kappa,1-\frac
4\kappa,\frac 8\kappa,r\right)$$
for $\kappa<8$. For $\kappa=6$, one recovers Cardy's formula. For
$\kappa\in (0,8/3]$, this can be interpreted as the probability that
two $\SLE_\kappa$'s connecting $x_1,x_2$ and $x_3,x_4$ resp. do not
intersect and are not connected by a loop in an independent loop-soup
with intensity $\lambda_\kappa$. 

A solution of the hypergeometric equation with
parameters $(a,b,c)=(4/\kappa,1-4/\kappa,8/\kappa)$ can be written as:
$$\int_{C'}
u^{-4/\kappa}(1-u)^{12/\kappa-2}(1-ux)^{-4/\kappa}du=x^{1-8/\kappa}\int_C
v^{-4/\kappa}(v-x)^{12/\kappa-2}(v-1)^{-4/\kappa}dv$$
where $C$ is a formal linear combination of ``cycles'', i.e. paths of
integrations over which the integrand is single-valued and are either
closed or start and end at singular points of the integrand, i.e.
($0,1,x,\infty$).
In terms of the original variables, one gets the following expression
(inverting the roles of $x_3$ and $x_4$):
\begin{equation}\label{en2}
\varphi(x_1,x_2,x_3,x_4)=\prod_{i<j<4}(x_j-x_i)^{2/\kappa}\prod_{i<4}(x_4-x_i)^{1-6/\kappa}\int_C\prod_{i<4}(v-x_i)^{-4/\kappa}(v-x_4)^{12/\kappa-2}dv  
\end{equation}
Note that although (\ref{sys}) is symmetric in the $(2n)$ variables
$x_1,\dots,x_{2n}$, the integrand is symmetric only in the first three
variables. In fact, one can interchange the role of $x_4$ and, say,
$x_3$ by a homographic change of variables. It is possible, but
apparently not very practical, to rewrite this integrand
in a completely symmetric fashion.
 
\subsection{Case $\kappa=2$}
If one specializes the above expression to $\kappa=2$, one gets a
rational function, since $_2F_1(2,-1,4,r)=1-r/2$; note that
$2r(1-r/2)=1-(1-r)^2$. The corresponding solution $\varphi$ is:
\begin{align*}
\varphi(x_1,x_2,x_3,x_4)&=\frac 1{(x_1-x_2)^2(x_3-x_4)^2}-\frac
1{(x_1-x_3)^2(x_2-x_4)^2}\\
&=\det\left(
\begin{array}{cc}
(x_1-x_2)^{-2} & (x_1-x_3)^{-2} \\ (x_4-x_2)^{-2} & (x_4-x_3)^{-2}
\end{array}\right)
\end{align*}
For a general value of $n$, it is possible (if a bit tedious) to check
that:
$$\varphi(x_1,\dots,x_{2n})=\det\left((x_i-x_{j+n})^{-2}\right)_{1\leq
i,j\leq n}$$
gives a solution of (\ref{sys}), corresponding to the pairing
$\{(x_1,x_{n+1}),\dots,(x_n,x_{2n})\}$. Other pairings gives different
solutions; there are non-trivial linear relations between those solutions.
These solutions correspond to the scaling limit of Fomin's
formulae (\cite{Fomin,KozL}). More precisely, multiplying the determinant by the product of its
diagonal terms, one gets an alternating sum of probabilities of a certain non-intersection
event. It is interesting to observe that if $(x_i,x_j)$ are paired,
the limit of this probability as $x_i\rightarrow x_j$ is given by a
lower dimensional determinant, as it should. Likewise, the
decorrelation of SLEs living on different scales is obvious from these
formulae. For the particular pairing
$(x_1,x_{2n}),\dots,(x_n,x_{n+1})$, if the points are in cyclic order,
the probability of the corresponding non-intersection event is given by a single determinant.

Each permutation of $x_1,\dots,x_{2n}$ in this formula give a solution
of the system; there are (many) non-trivial linear relations between those solutions.
In general, other (geometric) pairings do not correspond to a single
determinant, but to a linear combination of determinants of this
type. One has to determine the solution satisfying appropriate
boundary conditions. It is clear that as $x_i\rightarrow x_j$, each
determinant (divided by its diagonal terms) goes either to $0$ or to a
lesser order determinant not involving $x_i,x_j$.

\subsection{Case $\kappa\rightarrow\infty$}

Consider the following system:
$$\left\{\begin{array}{l}\displaystyle\partial_{kk}, \hspace{1cm} k=1,\dots, 2n\\
\sum_k\partial_k\\
\sum_k x_k\partial_k-n\\
\sum_k x_k^2\partial_k-(x_1+\cdots+ x_{2n})
\end{array}\right.
$$
which can be seen as the (somewhat degenerate) limit of (\ref{sys}) as $\kappa$ goes to
infinity. Then it is easy to see that solutions (functions annihilated
by those operators) are polynomials in
$x_1,\dots x_{2n}$ homogeneous
of degree $n$, of partial degree at most 1 in each variable. The
vector space of these polynomials has dimension $C(2n,n)$; the 
dimension of the subspace of such polynomials that are also translation invariant is
$$C(2n,n)-C(2n,n-1)=C(2n,n)/(n+1)=C_n.$$
Considering the derivative of
these polynomials with respect to a variable, and by induction on the
number of variables, one gets that this linear space is spanned by:
$$(x_1-x_{n+1})\cdots(x_n-x_{2n})$$
and polynomials obtained by action of the symmetric group on this one.
Note that these spanning polynomials also satisfy the last relation
in the system. So the solution space has dimension exactly $C_n$.

\subsection{Case $\kappa=6$, $2n=6$}

Here we consider the system satisfied by the crossing probabilities
for critical percolation in an hexagon with alternate boundary
conditions (or a simply connected domain with six points marked on the
boundary). In this situation, there are $C_3=5$ elementary crossing
events, with probabilities adding up to 1 (see Figure
\ref{hexconfig}).
A solution $\varphi$ of this system is conformally
invariant (and not only covariant); so one can send $x_1,x_2,x_6$ to
$\infty,0,1$ by an homography, and $\varphi$ becomes a function of
three real variables $y_1,y_2,y_3$. This solution must be annihilated by six
second-order operators; it turns out that there is (exactly) one linear relation
between these operators. Singularities for these equations occur when
two of the $y$ variables are equal, or one is equal to
$0,1,\infty$. Consider an open set $U$ of regular points. 
Let $D(U)$ be the subalgebra of
differential operators in
$\C(y_1,y_2,y_3)[\partial_1,\partial_2,\partial_3]$ with coefficients
regular in $U$, and  $I$ be the left ideal of differential
operators generated by these equations.

It is possible to write the system as an integrable, rank 5 Pfaffian
system (see e.g. \cite{Yosh}; this is also true for any $\kappa>0$). 
Also, note that the constant function
is a solution; we will use this simple fact to explicitate solutions of
the system. Consider the image of the the solution space in $U$
(one may also think in terms of local systems, sheaves of
algebras, ...). Its image under
the operator $\partial_1$ has dimension 4. Let $J$ be the left ideal:
$$J=\{L\in D(U): L\partial_1\in I\}.$$
Then $J$ is a holonomic ideal of rank 4. There is a well-known example
of a holonomic system of rank 4, with 3 variables, and the same
singular locus as $J$, namely the system of Lauricella's $F_D$
function, which is the natural generalization of Appell's $F_1$ (see
e.g. \cite{B,ApKamp}). 

Recall that the system $F_D$ with parameters
$a,\gamma,\beta_1,\dots \beta_m$ in the variables $u_1,\dots,u_m$, is
given by the linear operators:
$$\sum_{j=1}^m u_j(1-u_i)\partial_{ij}+(\gamma-(a+\beta_i+1)u_i)\partial_i-\sum_{j\neq
i}\beta_iu_j\partial_j-a\beta_i$$
for $i=1,\dots,m$. A solution then also satisfies the Euler-Darboux
equations (if $\gamma-a-1\neq 0$):
$$(u_i-u_j)\partial_{ij}-(\beta_j\partial_i-\beta_i\partial_j)$$
for $1\leq i<j\leq m$. It is known that this system is of rank $(m+1)$
and that an integral representation of solutions is given by:
$$\int_C t^{a-1}(1-t)^{\gamma-a-1}\prod_j(1-tu_j)^{-\beta_j}dt.$$

By elimination in $D(U)$, it is possible to
exhibit elements in $J$ that generate a rank 4 ideal, the solution
space of which is a conjugate of that of a particular $F_D$ system.
More precisely, with the help of a formal computation
software, one can prove that if $\varphi$ is a solution of the
original system, then:
$$\frac{(y_1(y_2-y_1)(y_3-y_1)(1-y_1))^{2/3}}{((1-y_2)(1-y_3)(y_3-y_2)y_2y_3)^{1/3}}\partial_1\varphi(y_1,y_2,y_3)$$
is a solution of the $F_D$ system with parameters:
$$a=1/3,\gamma=2/3,\beta_1=-4/3,\beta_2=2/3,\beta_3=2/3.$$
Similarly, in order to get a solution symmetric in the $y$ variables,
one can do the same computation for the operator
$(y_1\partial_1+y_2\partial_2+y_3\partial_3)$, which by conformal
invariance corresponds to a perturbation of $x_6$. Then it turns out
that:
$$\frac{((1-y_1)(1-y_2)(1-y_3))^{2/3}}{(y_1y_2y_3(y_2-y_1)(y_3-y_1)(y_3-y_2))^{1/3}}(y_1\partial_1+y_2\partial_2+y_3\partial_3)\varphi(y_1,y_2,y_3)$$
is a solution of the $F_D$ system with parameters:
$$a=1/3,\gamma=8/3,\beta_1=2/3,\beta_2=2/3,\beta_3=2/3.$$
Solutions of this system are given by
$$\int_C
u^{-2/3}(1-u)^{4/3}\left((1-uy_1)(1-uy_2)(1-uy_3)\right)^{-2/3}du$$
where as before $C$ is an adequate ``cycle''. To get to $\varphi$, we
write:
$$\frac{d}{ds}\varphi(sy_1,sy_2,sy_3)=((y_1\partial_1+y_2\partial_2+y_3\partial_3)\varphi)(sy_1,sy_2,sy_3)$$
which after trivial manipulations leads to the following expression
for $\varphi(x_1,\dots, x_6)$:
\begin{equation}\label{ek6}
\prod_{i<j<6}(x_j-x_i)^{1/3}\prod_{i<6}(x_6-x_i)^0\int_C\prod_{\substack{i\in\{1,2\}\\j<6}}
(u_i-x_j)^{-2/3}\prod_{i\in\{1,2\}}(u_i-x_6)^0(u_2-u_1)^{4/3}du_1du_2.
\end{equation}

We end this section by quoting a particular result for configurations
with symmetries. Consider the domains $(\U,1,u,j,ju,j^2,j^2u)$ where
$\U$ is the unit disk, $j=e^{2i\pi/3}$, and $u$ is on the arc
$(1,j)$. Restricting a solution $\varphi$ to such domains, one gets a
function $g$ of a single variable $u$. Now, computing in $D(U)$ as
above, it is possible to derive a third-order ODE satisfied by
$g$ (of course this ODE has no constant term). The fact that the rank goes down is linked to the particular
symmetries of the configurations we study, and is easily interpreted
in terms of percolation observables (different elementary events have
the same probability because of the threefold rotational symmetry of the configuration). Making the cubic change of
variables $v=u^3$ (that sends the singularities $1,j,j^2$ to $1$),
then the quadratic change $w=-(v-1)^2/4v$ (that sends the
singularities $0,\infty$ to $\infty$), 
so that $w=\sin^2(3\theta/2)$ if $u=e^{i\theta}$,
it turns out that $w\mapsto w^{1/2}g'(w)$ satisfies the (classical)
hypergeometric equation
with parameters $(a,b;c)=(5/6,5/6;7/6)$. Hence $g'$ belongs to
the vector space spanned by:
$$h_1(w)=w^{-1/2}\vphantom{F}_2F_1\left(\frac 56,\frac 56;\frac 32;1-w\right), 
h_2(w)=w^{-1/2}(1-w)^{-1/2}\vphantom{F}_2F_1\left(\frac 13,\frac 13;\frac 12;1-w\right)$$
Let $g_i(w)=\int_w^1h_i(s)ds$. The function $g$ belongs to the vector
space spanned by $({\bf 1},g_1,g_2)$. Let us work out boundary
conditions for crossing probabilities. Consider:
$$g(u)=\varphi((\U,1,u,\dots,j^2u))=\P\left((1,u)\leftrightarrow (j,ju) \leftrightarrow (j^2,j^2u)\right)$$
the probability that the three blue sides belong to the same blue
cluster (the ``Mercedes'' configuration). Let
$$g_+(w)=g(u)+g(j/u), g_-(w)=g(j/u)-g(u)$$
for $u$ in the arc $(1,e^{i\pi/3})$  (note that when $u$ goes from $1$ to $j$, $w$ goes from
$0$ to $1$ and then back to $0$).
For parity reasons, it is easy to see that one can write $g_+=c_1g_1+c_3$,
$g_-=c_2g_2$. We need to determine the three constants
$c_1,c_2,c_3$. They are fixed by the boundary conditions:
$g_+(0^+)=g_-(0^+)=1$, and $g_+(w)-g_-(w)=O(\sqrt w)=O(u-1)$ as
$w\searrow 0$ (the exponent here is three times the half-plane one-arm exponent). From the Euler integral:
$$_3F_2(a_1,a_2,a_3;\rho_1\rho_2;x)=B(a_1,\rho_1-a_1)^{-1}B(a_2,\rho_2-a_2)^{-1}\int_{[0,1]^2} s^{a_1-1}(1-s)^{\rho_1-a_1-1}
 t^{a_2-1}(1-t)^{\rho_2-a_2-1}(1-stx)^{-a_3}dsdt$$
one gets:
\begin{align*}
g_1(0)&=B(5/6,2/3)^{-1}\int_0^1 (1-s)^{-1/2}\int_0^1 t^{-1/6}(1-t)^{-1/3}(1-st)^{-5/6}dsdt\\
&=\vphantom{F}_3F_2(1,5/6,5/6;3/2,3/2;1)B(1,1/2)\\
g_2(0)&=B(1/3,1/6)^{-1}\int_0^1(s(1-s))^{-1/2}\int_0^1
t^{-2/3}(1-t)^{-5/6}(1-ts)^{-1/3}dsdt\\
&=\vphantom{F}_3F_2(1/2,1/3,1/3;1,1/2;1)B(1/2,1/2)=\vphantom{F}_2F_1(1/3,1/3;1;1)B(1/2,1/2)=\frac{\Gamma(1/3)\Gamma(1/2)^2}{\Gamma(2/3)^2}.
\end{align*}
For the condition $(c_1h_1-c_2h_2)(w)=O(w^{-1/2})$, we need the
analytic continuation formulae:
\begin{align*}
w^{1/2}h_1(w)&=\frac{\Gamma(3/2)\Gamma(-1/6)}{\Gamma(2/3)^2}\vphantom{F}_2F_1\left(\frac
56,\frac 56;\frac 76;w\right)+\frac{\Gamma(3/2)\Gamma(1/6)}{\Gamma(5/6)^2}w^{-1/6}\vphantom{F}_2F_1\left(\frac
23,\frac 23;\frac 56;w\right)\\
w^{1/2}h_2(w)&=\frac{\Gamma(1/2)\Gamma(-1/6)}{\Gamma(1/6)^2}\vphantom{F}_2F_1\left(\frac
56,\frac 56;\frac 76;w\right)+\frac{\Gamma(1/2)\Gamma(1/6)}{\Gamma(1/3)^2}w^{-1/6}\vphantom{F}_2F_1\left(\frac
23,\frac 23;\frac 56;w\right).\\
\end{align*}
Hence we can conclude that:
$$c_2=\frac{\Gamma(2/3)^2}{\Gamma(1/3)\Gamma(1/2)^2}=\frac{\sqrt
3}{2\pi^2}\Gamma(2/3)^3,
c_1=c_2\frac{2\Gamma(5/6)^2}{\Gamma(1/3)^2}=\left(\frac{\sqrt 3}{2^{2/3}\pi}\right)^5\Gamma(2/3)^9.$$

\begin{figure}[htbp]
\begin{center}
\centerline{\psfig{file=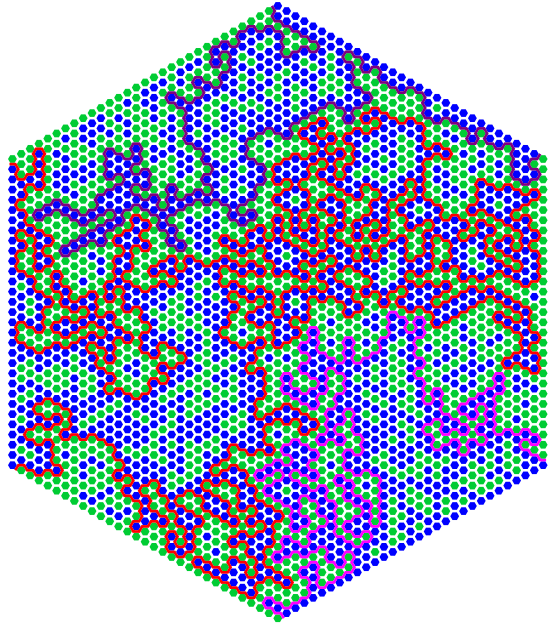}}
\end{center}
\caption{Percolation in a regular hexagon}\label{reghex}
\end{figure}

So the probability of a configuration where two (given) blue sides and two
yellow sides are connected for a regular
hexagon ($u=e^{i\pi/6}$, $w=1$, see Figure \ref{reghex}) is:
$$\frac{1-c_3}3=\frac{c_1g_1(0)}3=\frac 23\left(\frac{\sqrt 3}{2^{2/3}\pi}\right)^5\Gamma(2/3)^9\vphantom{F}_3F_2(1,5/6,5/6;3/2,3/2;1).$$

%$f(w)=g(w)w^{1/2}$ satisfies
%the ODE:
%$$72w^2(1-w)f'''(w)+w(264-372w)f''(w)+(126-338w)f'(w)-25f(w)=0$$
%which, up to a multiplicative constant, is the ODE associated with the
%hypergeometric function:
%$$_3F_2\left(\frac 12,\frac 56,\frac 56;\frac 76,\frac 32;w\right)$$
%Note that the regular hexagon corresponds to $w=1$, so that there
%may be a connection between values of percolation probabilities for a
%regular hexagon and evaluation of some hypergeometric functions at unity
%(which is the case for Watts' formula in the square; see \cite{Maier}).

\section{Formal solutions}

In this section, we discuss a family of formal solutions for the
system (\ref{sys}). These are integrals taken on a ``cycle''; they
define actual solutions under conditions of convergence and
determination of the integrand. We shall
discuss later the problem of identifying probabilistic solutions.

Consider the function:
\begin{align*}
\phi({\bf x},{\bf u})=&\prod_{\substack{1\leq i< n\\1\leq j<2n}}
\left(u_i-x_j\right)^{-\frac 4\kappa}\prod_{1\leq
i<n}\left(u_i-x_{2n}\right)^{\frac{12}\kappa-2}\prod_{1\leq
i_1<i_2<n}\left(u_{i_2}-u_{i_1}\right)^{\frac 8\kappa}\\
&\prod_{1\leq
j_1<j_2<2n}\left(x_{j_2}-x_{j_1}\right)^{\frac 2\kappa}\prod_{1\leq
j<2n}\left(x_{2n}-x_{j}\right)^{1-\frac 6\kappa}
\end{align*}
where ${\bf x}=(x_1,\dots,x_{2n})$, ${\bf u}= (u_1,\dots u_{n-1})$.
 We want
to prove that a function of the parameters $(x_1,\dots x_{2n})$ of the form:
$$\int_C\phi({\bf x},{\bf u})d{\bf u}$$
defines a solution of (\ref{sys}) for appropriate cycles $C$ (that may
depend on the $x$ parameters). Note that this covers the cases
(\ref{en2}) and (\ref{ek6}). The conformal covariance conditions are
easy to check. Denote by ${\mc L}_k$ the differential operator:
$${\mc L}_k=\frac\kappa 2\partial_{kk}+\sum_{l\neq
k}\frac{2\partial_{l}}{x_l-x_k}+\frac{\kappa-6}\kappa\sum_{l\neq
k}\frac 1{(x_l-x_k)^2}.$$  
Then the following lemma holds.

\begin{Lem}
(i) If $1\leq k<2n$, then:
$${\mc L}_k\phi=-\sum_{1\leq i<n}\frac{\partial}{\partial u_i}\left(\frac{2\phi}{u_i-x_k}\right)$$
(ii) Moreover:
$${\mc L}_{2n}\phi=-\sum_{1\leq i<n}\frac{\partial}{\partial
u_i}\left(\left(\frac
2{u_i-x_{2n}}+\frac{\kappa-8}{u_i-x_{2n}}\prod_{1\leq
j<2n}\frac{u_i-x_j}{x_{2n}-x_j}\prod_{1\leq j<n,j\neq i}\left(\frac{x_{2n}-u_j}{u_i-u_j}\right)^2\right)\phi\right)$$
\end{Lem}
\begin{proof}
(i) By symmetry, one can assume that $k=1$. Then:
\begin{align*}
\frac{{\mc L}_1\phi}\phi=&
\frac\kappa 2\left(\left(-\frac 4\kappa\sum_{1\leq i<n}\frac
1{x_1-u_i}+\frac 2\kappa\sum_{1<j<2n}\frac 1{x_1-x_j}+\left(1-\frac
6\kappa\right)\frac 1{x_1-x_{2n}}\right)^2\right.\\
&\left. +\frac 4\kappa\sum_{1\leq i<n}\frac
1{(x_1-u_i)^2}-\frac 2\kappa\sum_{1<j<2n}\frac 1{(x_1-x_j)^2}-\left(1-\frac
6\kappa\right)\frac 1{(x_1-x_{2n})^2}\right)\\
&+\sum_{1<j<2n}\frac2{x_j-x_1}\left(-\frac 4\kappa\sum_{1\leq i<n}\frac
1{x_j-u_i}+\frac 2\kappa\sum_{\substack{1\leq k<2n\\k\neq j}}\frac 1{x_j-x_k}+\left(1-\frac
6\kappa\right)\frac 1{x_j-x_{2n}}\right)\\
&+\frac 2{x_{2n}-x_1}\left(1-\frac 6\kappa\right)\left(-2\sum_{1<i<n}\frac
1{x_{2n}-u_i}+\sum_{1\leq j<2n}\frac
1{x_{2n}-x_j}\right)+\left(1-\frac 6\kappa\right)\sum_{1<j\leq
2n}\frac 1{(x_j-x_1)^2}
\end{align*}
which after simplifications leads to:
\begin{align*}
\frac{{\mc L}_1\phi}\phi=&
-2\left(1+\frac 4\kappa\right)\sum_{1\leq i<n}\frac
1{(x_1-u_i)^2}-\frac 8\kappa\sum_{\substack{1\leq
i<n\\1<j<2n}}\left(\frac 1{x_1-u_i}\frac 1{x_1-x_j}+\frac
1{x_j-x_1}\frac 1{x_j-u_i}\right)\\
&-4\left(1-\frac 6\kappa\right)\sum_{1\leq i< n}\left(\frac 1{x_1-x_{2n}}\frac 1{x_1-u_i}+\frac
1{x_{2n}-x_i}\frac 1{x_{2n}-u_i}\right)\\
&+\frac {16}\kappa\sum_{1\leq i_1<i_2<n}\frac
1{x_1-u_{i_1}}\frac 1{x_1-u_{i_2}}
\end{align*}
where we have used the identity:
$$\frac 1{(a-b)(a-c)}+\frac 1{(b-a)(b-c)}+\frac 1{(c-a)(c-b)}=0$$
for $(a,b,c)=(x_1,x_i,x_j)$, $1<i<j\leq 2n$. Also:
\begin{align*}
\frac 1\phi\frac{\partial}{\partial u_i}\left(\frac{2\phi}{u_i-x_1}\right)
&=\frac 2{u_i-x_1}\left(-\left(1+\frac 4\kappa\right)\frac 1{u_i-x_1}-\frac
4\kappa\sum_{1<j<2n}\frac 1{u_i-x_j}\right.\\
&\left.+\left(\frac
{12}\kappa-2\right)\frac 1{u_i-x_{2n}}+\frac
8\kappa\sum_{\substack{1\leq i_2<n\\i_2\neq i}}\frac 1{u_i-u_{i_2}}\right)
\end{align*}
Making use of the same identity as above, one can easily conclude.

(ii) The situation here is slightly more intricate than in
(i). Let $P$ be the rational function:
$$P=\frac 1\phi\left({\mc L}_{2n}\phi+\sum_{1\leq i<n}\frac{\partial}{\partial
u_i}\left(\frac{2\phi}{u_i-x_{2n}}\right)\right)$$
Reasoning as above, one gets the following identity:
\begin{align*}
P=&
\frac{(8-\kappa)(4-\kappa)}\kappa\left(\sum_{1\leq i<n}\frac
3{(u_i-x_{2n})^2}
-2\sum_{\substack{1\leq i<n\\1\leq j<2n}}\frac{1}{(u_i-x_{2n})(x_j-x_{2n})}
\right.\\&\left.
+\sum_{1\leq j
_1<j_2<2n}\frac 1{(x_{j_1}-x_{2n})(x_{j_2}-x_{2n})}+4\sum_{1\leq i_1<i_2<n}\frac 1{(u_{i_1}-x_{2n})(u_{i_2}-x_{2n})}
\right).
\end{align*}
On the other hand, if $Q$ is the rational function:
$$Q=\frac 1\phi\sum_{1\leq i<n}\frac{\partial}{\partial
u_i}\left(\left(\frac{8-\kappa}{u_i-x_{2n}}\prod_{1\leq
j<2n}\frac{u_i-x_j}{x_{2n}-x_j}\prod_{1\leq j<n,j\neq i}\left(\frac{x_{2n}-u_j}{u_i-u_j}\right)^2\right)\phi\right)$$
then $Q$ can be written as:
\begin{align*}
Q=&
\frac{(8-\kappa)(4-\kappa)}\kappa
\sum_{1\leq i<n}\frac{1}{u_i-x_{2n}}\prod_{1\leq
j<2n}\frac{u_i-x_j}{x_{2n}-x_j}
\prod_{1\leq j<n,j\neq i}\left(\frac{x_{2n}-u_j}{u_i-u_j}\right)^2\\
&\times\left(-\sum_{1\leq j<2n}\frac 1{u_i-x_j}+\frac 3{u_i-x_{2n}}+
2\sum_{1\leq j<n,j\neq i}\frac 1{u_i-u_j} 
\right).
\end{align*}
So the statement reduces to the identity of rational functions
$P=Q$. This is a bit tedious, and we include these computations for
the sake of completeness.
We can assume that $\kappa\neq 4,8$. Note also that $P$ and $Q$ are symmetric in the $u$ variables. Let
us expand $Q$ at $u_1=x_{2n}$ (denote $\eps=u_1-x_{2n}$):
\begin{eqnarray*}
Q&=&\frac{(8-\kappa)(4-\kappa)}\kappa
\frac 1{\eps}\left(1+\eps\sum_{1\leq j<2n}\frac
1{x_{2n}-x_j}+\eps^2\sum_{1\leq j_1<j_2<2n}\frac
1{(x_{2n}-x_{j_1})(x_{2n}-x_{j_2})}\right)\\
&&\times\left(1+\eps\sum_{1<i<n}\frac {-2}{x_{2n}-u_i}+\eps^2\left(\sum_{1<i_1<i_2<n}\frac
4{(x_{2n}-u_{i_1})(x_{2n}-u_{i_2})}+\sum_{1<i<n}\frac 3{(x_{2n}-u_i)^2}\right)\right)\\
&&\times\left(\frac 3\eps+2\sum_{1<i<n}\frac
1{x_{2n}-u_i}-\sum_{1\leq j<2n}\frac 1{x_{2n}-x_j}\right.\\
&&\hphantom{\times}\left.+\eps\left(
\sum_{1\leq j<2n}\frac 1{(x_{2n}-x_j)^2}-2\sum_{1<i<n}\frac 1{(x_{2n}-u_i)^2}\right)\right)+o(\eps)\\
&=&\frac 3{\eps^2}+\frac 1\eps\left(\sum_{1\leq
j<2n}\frac{3-1}{x_{2n}-x_j}+\sum_{1<i<n}\frac
{2-6}{x_{2n}-u_i}\right)\\
&& +\sum_{1\leq
j_1<j_2<2n}\frac{3-2}{(x_{2n}-x_{j_1})(x_{2n}-x_{j_2})}+\sum_{1<i_1<i_2<n}\frac{12-8}{(x_{2n}-u_{i_1})(x_{2n}-u_{i_2})}\\
&& +\sum_{1<i<n}\frac{9-2-4}{(x_{2n}-u_i)^2}+\sum_{1<j<2n}\frac{1-1}{(x_{2n}-x_j)^2}+\sum_{\substack{1<i<n\\1\leq
j<2n}}\frac{-2+2-2}{(x_{2n}-x_j)(x_{2n}-u_i)}+o(\eps).
\end{eqnarray*}
From here, and using symmetry, it appears readily that: 
$$P-Q=o((u_i-x_{2n}))$$
for $1\leq i<n$. One can also check that $Q$ is regular along $u_i=u_j$,
for $1\leq i<j<n$ (by symmetry it is enough to check that the term in
$(u_i-u_j)^{-2}$ vanishes). Considering $(P-Q)$ as a rational function of the
$u$ variables, it appears that it has no pole (even at infinity), so
it is constant (i.e. a function only of the $x$ variables). Setting
$u_1=x_{2n}$, it appears that $(P-Q)$ is identically zero, which
concludes the proof.
\end{proof}

One can rephrase the lemma as follows: there are differential $(n-2)$-forms
$\omega_1,\dots ,\omega_{2n}$ with rational coefficients such that for
$1\leq k\leq 2n$:
$${\mc L}_k (\phi du_1\wedge\dots\wedge du_{n-1})=d(\phi\omega_k).$$
So for a ``cycle'' $C$, we get a (real) solution of (\ref{sys}) as
soon as the following formal computation makes sense:
$${\mc L}_k\int_C \phi du_1\wedge\dots\wedge du_{n-1}=\int_C {\mc
L}_k( \phi du_1\wedge\dots\wedge du_{n-1})=\int_C
d(\phi\omega_k)=\int_{\partial C}\phi\omega_k=0.$$
The cycle $C$ must be such that $\phi$ has a single-valued
determination on $C$, one can differentiate w.r.t. the $x$ parameters in the integral, 
and $\phi\omega_k$ vanishes (to a sufficient order) on $\partial
C$ (which is true in particular if $\partial C=\varnothing$). In fact $C$ is a cycle for the twisted de Rham homology associated
with $\phi$ (see e.g. \cite{SST}, section 5.4).
In the next section we give examples of such cycles, in relation with
the probabilistic set-up.

\section{Explicit solutions}

In this section we propose choices of cycles of integration such that
one gets well-defined solutions, and that these solutions can be
interpreted in the probabilistic situations discussed earlier. We
distinguish the cases $\kappa\in (0,8/3]$ and $\kappa=6$. Recall that
(with an emphasis on the number of parameters):
\begin{align*}
\phi_n({\bf x},{\bf u})=&\prod_{\substack{1\leq i< n\\1\leq j<2n}}
\left(u_i-x_j\right)^{-\frac 4\kappa}\prod_{1\leq
i<n}\left(u_i-x_{2n}\right)^{\frac{12}\kappa-2}\prod_{1\leq
i_1<i_2<n}\left(u_{i_2}-u_{i_1}\right)^{\frac 8\kappa}\\
&\prod_{1\leq
j_1<j_2<2n}\left(x_{j_2}-x_{j_1}\right)^{\frac 2\kappa}\prod_{1\leq
j<2n}\left(x_{2n}-x_{j}\right)^{1-\frac 6\kappa}
\end{align*}
where ${\bf x}=(x_1,\dots,x_{2n})$, ${\bf u}=(u_1,\dots u_{n-1})$.

\begin{figure}[htbp]
\begin{center}
\scalebox{.5}{\input{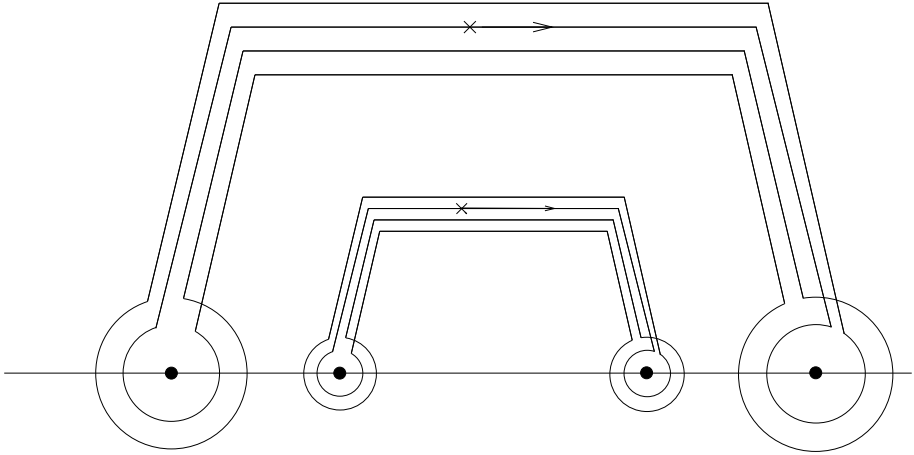}}
\end{center}
\caption{Two non-intersecting double contour loops}\label{doublecont}
\end{figure}

Assume that $x_1<\cdots<x_{2n}$ are boundary points of $\H$, and the
involution $\iota$ with no fixed points define a non-crossing pairing
of the $(2n)$ points $x_1\dots x_{2n}$. Recall that the fundamental
group $\Pi_1(\C\setminus\{x_k,x_{\iota(k)}\},z_0)$ is the free group generated by two
elements $\sigma_k,\sigma_{\iota(k)}$ corresponding to loops around
$x_k,x_{\iota(k)}$ respectively. A double contour loop around
$x_k,x_{\iota(k)}$ corresponds to the commutator
$$\sigma_k\sigma_{\iota(k)}\sigma_k^{-1}\sigma_{\iota(k)}^{-1}$$
so that its image in $\Pi_1(\C\setminus\{x_k\},z_0)$ and
$\Pi_1(\C\setminus\{x_{\iota(k)}\},z_0)$ is the identity.
We assume that the other $x$ points are outside this loop. Then $\phi_n$
is single-valued on such a double contour loop. 
Now let us enumerate the $(n-1)$ pairs that do not contain the last
point $x_{2n}$:
$$\{\{x_k,x_{\iota(k)}\},k=1,\dots,2n\}=\{\{a_k,b_k\},k=1,\dots (n-1)\}\cup\{\{x_{2n},x_{\iota(2n)}\}\}.$$ 
Denote also $a_n=x_{\iota(2n)}$, $b_n=x_{2n}$, and assume that
$a_k<b_k$, $k=1,\dots n$.
For $k=1,\dots,n-1$, let $C_k$ be a double contour loop around
$a_k,b_k$, so that these loops do not intersect pairwise (see Figure
\ref{doublecont}). 
Define $C$ to be the Cartesian product of these
loops:
$$C=C_1\times\cdots\times C_{n-1}.$$
Of course, $C$ is a function of $x_1,\dots,x_{2n}$. Then $\phi_n$ (as a function of $u_1,\dots,u_{n-1}$) is single-valued
on $C$, and $\partial C=\varnothing$. Moreover, the integral of
$\phi_n$ on $C$ does not change if the loops are deformed.

If $L$ is a double contour loop around $0,1$, and $p,q$ are numbers
with real part larger than $(-1)$, then:
\begin{align*}
\int_Lt^p(1-t)^qdt=&\left(1-e^{2i\pi q}+e^{2i\pi (p+q)}-e^{2i\pi
p}\right)\int_0^1t^p(1-t)^qdt\\
&=(1-e^{2i\pi p})(1-e^{2i\pi q})\frac{\Gamma(p+1)\Gamma(q+1)}{\Gamma(p+q+2)}
\end{align*}
as is easily seen when $L$ gets close to the unit segment. By analytic
continuation in $p$ and $q$, the second expression stays valid for
general values of $p$ and $q$. In the situation where $p=q$, one can
use figure eight loops.

Assume that the two consecutive points $x_k,x_{k+1}$, $k+1<2n$, are paired and
correspond to the loop $C_j$. Then it is easy to see that:
$$\lim_{x_{k+1}\rightarrow x_k}(x_{k+1}-x_k)^{6/\kappa-1}\int_C\phi_n({\bf x},{\bf u}) du_1\dots
du_{n-1}=c\int_{\widehat C}\phi_{n-1}(\hat{\bf x},\hat{\bf u})du_1\dots \widehat {du_j}\dots du_{n-1}$$ 
where $\hat{\bf x}=(x_1,\dots,\widehat{x_k},\widehat{x_{k+1}},\dots,x_{2n})$,
$\hat{\bf u}=(u_1,\dots,\widehat{u_j},\dots,u_{n-1})$, and
$\widehat C=C_1\times\cdots\widehat{C_j}\times\cdots\times C_{n-1}$. The constant
$c$ is given by (up to a complex number of modulus one depending on
the choice of a determination for $\phi_n$ on $C$):
$$c=c_\kappa=4\sin\left(\frac
4\kappa\pi\right)^2\frac{\Gamma\left(1-\frac
4\kappa\right)^2}{\Gamma\left(2-\frac 8\kappa\right)}=\frac
{4\pi^2}{\Gamma\left(2-\frac 8\kappa\right)\Gamma\left(\frac 4\kappa\right)^2}$$
This constant is non zero, since we assume that $(8/\kappa)\notin\N$.
If $n=2$, then the limit is $c(x_4-x_3)^{1-6/\kappa}$.

In this set-up, we can now formulate:

\begin{Prop}
Let $\kappa\in(0,8/3)$, $8/\kappa\notin\N$. Let
$\gamma_1,\dots,\gamma_n$ be (the traces of) $n$ independent
$\SLE_\kappa$'s from $a_k$ to $b_k$, $k=1,\dots n$, and $L_2,\dots
L_{n}$ be independent loop-soups with intensity $\lambda_\kappa$. Then:
$$\P((\gamma_j)^{L_j}\cap
(\cup_{i<j}\gamma_i=\varnothing),j=2,\dots,n)=(c_\kappa)^{1-n}\prod_{k=1}^n(b_k-a_k)^{6/\kappa-1}\left|\int_C\phi_n({\bf
x},{\bf u})d{\bf u}\right|$$
\end{Prop}
\begin{proof}
Denote by $\tilde\psi(x_1,\dots x_{2n})$ the right-hand side of the equation
in the proposition. Then we want to prove $\tilde\psi=\psi$ (with the
notations of Proposition \ref{Prestr}). By construction, $\tilde\psi$
has the property (ii) of Proposition \ref{Prestr}. Also, it is easy to
see that $\tilde\psi$ is conformally invariant. So let us momentarily
fix the values $x_1=0,x_{2n-1}=1,x_{2n}=\infty$. Then
$\tilde\psi$ extends continuously to the compactification of
$$\{(x_2,\dots,x_{2n-2})\in\R^{2n-3}\textrm{ s.t. } 0=x_1<\cdots< x_{2n-1}=1\}$$
so it is bounded. There is a $k\in\{1,\dots,n-1\}$ such that $a_k$ and
$b_k$ are two consecutive points on the real line. Consider the
martingale associated with $\tilde\psi$ and the evolution of the
$k$-th SLE. We use the hydrodynamic normalization at infinity;
$(g_t)$ are the conformal equivalences associated with the SLE, which
is defined up to time $\tau=\tau(b_k)$ (which is finite due to this
arbitrary choice of normalization). Then $g_\tau$ defines a
conformal equivalence between the unbounded connected component of
$\H\setminus \gamma_k$ and $\H$. As
$t\nearrow\tau$, $W_t$ and $g_t(b_k)$ go to the same finite limit, and
$g_t(a_j),g_t(b_j)$ have finite (distinct) limits. So the martingale goes to:
$$\tilde\psi_{n-1}(g_\tau(\hat{\bf x}))\prod_{j\neq k}\left(g_\tau'(a_j)g_\tau'(b_j)\left(\frac{b_j-a_j}{g_\tau(b_j)-g_\tau(a_j)}\right)^2\right)^{\alpha_\kappa}$$ 
with the obvious notations. The product can be written as:
$$\prod_{j\neq k}\P(\gamma_j^L\cap\gamma_k\neq\varnothing)$$
where $L$ is an independent loop-soup.
Now, as we remarked earlier, $\psi({\bf x})$
does not depend on the ordering of
the SLEs (this follows from the Poissonian nature of the loop
soup). So we can assume that the $(n-1)$ SLEs loop-soups are reordered so
that the last one corresponds to $k$ (that is, $k$ is re-indexed as $n$). 
As in Proposition \ref{Prestr}, the restriction property for the
loop-soup and for $(n-1)$ SLEs now implies that:
$$\psi({\bf x})=\E\left(\P\left((\tilde\gamma_j)^{\tilde L_j}\cap
(\cup_{i<j}\tilde\gamma_i)=\varnothing,j=2\dots (n-2)|g_\tau(a_j),g_\tau(b_j),j<n\right)\prod_{j<n}\ind_{(\gamma_j)^{L'_j}\cap\gamma_n=\varnothing}\right)$$
where $\tilde\gamma_j$ is an SLE from $g_\tau(a_j)$ to $g_\tau(b_j)$,
and $\tilde L_j$, $L'_j$ are independent loop-soups. Comparing with the
limiting value of the martingale as $t\nearrow\tau$, one concludes by
induction on $n$, using the optional stopping theorem.
\end{proof}

Note that it appears to be difficult to determine analytically the
boundary behaviour of $\tilde\psi$ along all boundary
components. The point here is that the dynamics of the SLEs lead to an exit
on a fixed boundary component; and commutation allows us to consider
the SLEs in an appropriate order. 

We gave probabilistic interpretations in the case $\kappa\leq 8/3$,
$\kappa=6$, $\kappa=8$. We have identified corresponding integration
cycles when $\kappa\leq 8/3$, $8/\kappa\neq\N$ and $\kappa=8$. Let us
discuss the remaining cases $\kappa=8/m$, $m\geq 3$, and $\kappa=6$.

The problem is to identify a cycle of integration corresponding to a
geometric configuration of SLEs. In particular, one can consider
the case of $n$ nested paths connecting $x_i$ to $x_{2n+1-i}$,
$i=1,\dots,n$, where $x_1,\dots, x_{2n}$ are in cyclical order on the
boundary. This is the natural situation for Fomin's formulae (but not
for the UST e.g.). Then the following cycle can be used (say when $\kappa<4$):
$$C=C_1\times\cdots\times C_{n-1}$$
where $C_i$ is a loop starting and ending at $x_{2n}$, circling
counterclockwise around $x_{2n-i}$, and leaving the other marked
points on its right-hand side; moreover, the $C_i$'s are chosen so that
they do not intersect except at $x_{2n}$, consequently the integrand
has a single-valued determination on $C$.

In the case where $\kappa=2$, induction on $n$ and a residue
computation give the determinant obtained as the scaling limit of
Fomin's formula. One could also have chosen $C'=C'_1\times\cdots\times
C'_{n-1}$, with $C_i$ a loop starting at $x_1$ and circling around
$x_{1+i}$ (and the special point in the integrand is now $x_1$, not
$x_{2n}$); the solution is then proportional to the one corresponding
to $C$. Here, using induction on $n$ (and analytic continuation in
$\kappa$), we can show that these two solutions are indeed
proportional and have the correct asymptotic behaviour when two
consecutive marked points collapse. Finding a cycle with geometrically
prescribed boundary conditions (in the Weyl chamber
$\{x_1<\cdots<x_{2n}\}$) seems to be technically difficult in general.

Also, we remark here that the integral representation becomes somewhat
degenerate in the particularly interesting case $\kappa=6$. As in the
case of standard hypergeometric functions, analytic continuation in
the parameters (here $\kappa$) can be put to good use. Since the
functions $\psi$ can in this instance be interpreted as (scaling
limits of) crossing probabilities, it would be nice to find real
integration cycles (formal sums of Cartesian products of segments
$[x_i,x_{i+1}]$) corresponding to these events (i.e. satisfying the
right boundary conditions; in general, when two consecutive points
collapse, the boundary condition is either $0$ or a crossing
probability for a configuration with $2(n-1)$ marked points). For
small values of $n$ (e.g. $2n=6$), one can work out such cycles,
though the general pattern is not very clear. 
Again for small $n$, these cycles give
enough solutions, so the results of this paper can be seen as giving
an algorithmic solution to the problem of computing crossing
probabilities for the scaling limit of percolation with alternate
boundary conditions.

As for general properties of the holonomic system we studied, the main
question is to prove that the rank is indeed $C_n$. From the case
$\kappa\rightarrow\infty$, we see that the rank should be at most
$C_n$. Under smoothness assumptions, the restriction construction give
$C_n$ distinct solutions when $0<\kappa\leq 8/3$; this should be
enough to prove that the rank is $C_n$ in general (since all these
systems then have a rank $C_n$ Pfaffian form that should have an
analytic continuation in $\kappa$).

\noindent{\bf Acknowledgments.} I wish to thank David Wilson and Greg
Lawler for fruitful conversations, as well as anonymous referees for
helpful comments.

\bibliographystyle{abbrv}
\bibliography{biblio}

-----------------------

Courant Institute of Mathematical Sciences

251 Mercer St., New York NY 10012

dubedat@cims.nyu.edu

\end{document}